\documentclass[11pt,letterpaper]{article}

\usepackage{amssymb,amsmath,amsfonts}
\usepackage{graphicx,xcolor,enumitem}
\usepackage{epsfig}
\usepackage{amsthm}
\usepackage{bm}
\usepackage{subcaption}
\usepackage{nicematrix}
\usepackage[numbers]{natbib}
\usepackage{csquotes}
\renewcommand{\mkbegdispquote}[2]{\itshape}
\usepackage[letterpaper, twoside, hmarginratio=1:1, vmarginratio=1:1, left=1in,top=1in]{geometry}

\RequirePackage[breaklinks=true, hidelinks]{hyperref}
\usepackage{breakcites}

\usepackage{algorithm,algorithmicx}
\usepackage{algpseudocode}

\usepackage{accents}

\newcommand{\cB}{\mathcal{B}}
\newcommand{\E}{\mathbb{E}}
\newcommand{\R}{\mathbb{R}}
\newcommand{\p}{\mathbb{P}}
\newcommand{\barp}{\bar{\mathbb{P}}}

\newcommand{\bS}{\mathbb{S}}
\newcommand{\cL}{{\mathcal L}}

\newcommand{\tr}{{\rm{tr}}}

\newcommand{\cV}{{\mathcal V}}
\newcommand{\cW}{{\mathcal W}}

\newtheorem{theorem}{Theorem}

\newtheorem{lemma}[theorem]{Lemma}

\theoremstyle{definition}

\numberwithin{equation}{section}
\numberwithin{theorem}{section}

\begin{document}

\title{Distributionally robust Kalman filtering with volatility uncertainty}

\author{Bingyan Han \thanks{Email: \href{mailto:bingyanhan@hkust-gz.edu.cn}{bingyanhan@hkust-gz.edu.cn}. Thrust of Financial Technology, The Hong Kong University of Science and Technology (Guangzhou), China. This work began when the author was a postdoctoral researcher in the Department of Mathematics at the University of Michigan. He expresses gratitude to the University of Michigan for providing support and an atmosphere conducive to this work. The author is also partially supported by the HKUST (GZ) Start-up Fund G0101000197, the Guangzhou-HKUST(GZ) Joint Funding Program (No. 2024A03J0630), and the National Natural Science Foundation of China (Grant No. 12401621). The author thanks the anonymous referees and editors for their valuable comments that have greatly improved this manuscript. This arXiv version has more calculation details than the published version, due to page limits.}
}

\date{May 30, 2025}
\maketitle

\begin{abstract}
	This work presents a distributionally robust Kalman filter to address uncertainties in noise covariance matrices and predicted covariance estimates. We adopt a distributionally robust formulation using bicausal optimal transport to characterize a set of plausible alternative models. The optimization problem is transformed into a convex nonlinear semi-definite programming problem and solved using the trust-region interior point method with the aid of $LDL^\top$ decomposition. The empirical outperformance is demonstrated through target tracking and pairs trading.
	\\[2ex] 
	\noindent{\textbf {Keywords}: Bicausal optimal transport, robust Kalman filtering, minimax optimization.}
	\\[2ex]
	%\noindent{\textbf {Mathematics Subject Classification:} 91G80, 91A80, 60G22, 60H20.} \\
	%\noindent{\textbf {JEL Classification:} C72, C73, G11.} 
\end{abstract}

\section{Introduction}
The Kalman filter is a widely used algorithm for estimating unknown state variables from noisy observations. It finds applications in diverse fields such as control systems, communications, statistics, and finance. While the standard Kalman filter assumes a linear system with Gaussian noise, this assumption is often unrealistic and the estimation of noise covariance matrices is challenging. Consequently, the classical Kalman filter is prone to model misspecification, resulting in a substantial body of literature on robust Kalman filters. Various robust formulations exist, including $H^\infty$ filtering \cite{petersen1999robust,hassibi1999indefinite}, risk-sensitive filtering \cite{speyer1974optimization}, outlier-robust approaches \cite{durovic1999robust,gandhi2009robust,huang2020novel}, and distributionally robust methods utilizing the Kullback-Leibler (KL) divergence \cite{levy2012robust,zorzi2016robust} or the optimal transport (OT) framework \cite{nips2018kalman}. Another relevant area is adaptive Kalman filtering; recent review by \cite{zhang2020identification} provides insights into this topic.

This paper tackles the problem of inaccurate or uncertain covariance matrices in Kalman filtering using the distributionally robust method. Since the KL divergence $D_{KL}(\p || \barp)$ is not symmetric and requires $\p$ to be absolutely continuous with respect to $\barp$, we opt for the Wasserstein distance instead. It measures the distance between two distributions, quantifying how difficult it is to transform a source probability into a target one. The OT framework offers a clear geometric interpretation and has found applications across diverse fields such as mathematics, machine learning, statistics, control theory, and economics; see \cite{villani2009optimal,taghvaei2020optimal,zorzi2020optimal} for an incomplete list.

Our method differs from previous research \cite{huang2017novel,nips2018kalman} because it utilizes a recent OT framework with (bi)causality constraints. When handling temporal data, OT needs to be formulated to reflect the sequential nature of the data. \cite{lassalle2013causal,backhoff2017causal} introduced Causal Optimal Transport (COT), imposing a causality constraint on the transport plans between discrete stochastic processes. Informally, given the past of a process $Z$, the past of another process $Z'$ should be independent of the future of $Z$ under the transport plan. If this condition also holds when the positions of $Z$ and $Z'$ are exchanged, it is called bicausal optimal transport (BCOT). 

BCOT is integrated into robust Kalman filtering as follows: At a given time $t$, uncertainty arises from inaccurate information about the covariance matrices of noises $Q_t$ and $R_t$. Additionally, the previous estimate $\Sigma_{t-1}$ of the covariance for the unobserved state $x_{t-1}$ given measurements $y_{1:t-1}$ can also be misspecified. We seek to devise a robust filter capable of rectifying potential errors in the previous estimate $\Sigma_{t-1}$ when the current measurement $y_t$ becomes available. Since two time steps, $t-1$ and $t$, are involved, it is crucial to appreciate the evolution of information flows by imposing the bicausality condition.

In our BCOT framework, alternative models are assumed to have a linear Gaussian structure, potentially with different covariance matrices $\bar{Q}_t$, $\bar{R}_t$, and $\bar{\Sigma}_{t-1}$. We derive an explicit formula for our BCOT framework using the dynamic programming principle \cite{backhoff2017causal}. For Gaussian distributions, an explicit formula for the Wasserstein distance with the $L_2$-norm exists in the one-step case \cite{bures1969extension,givens1984class}. The minimax formulation reduces to a convex nonlinear semi-definite programming problem, which can be solved numerically in various ways. We use the $LDL^\top$ decomposition method from \cite{benson2003solving} to reformulate the semi-definite constraints as more tractable convex constraints.

We demonstrate the effectiveness of our method through target tracking and pairs trading. Our framework accommodates any Kalman-like filter as the reference model. We employ the expectation--maximization (EM) algorithm to calibrate the reference model parameters and validate them with out-of-sample tests. Our results show that the BCOT method consistently outperforms the KL method \cite{zorzi2016robust}, the classical OT method \cite{nips2018kalman}, and the non-robust counterpart, providing strong empirical evidence for the benefit of the BCOT approach.

The paper is organized as follows. Section \ref{sec:form} reviews the classical Kalman filtering procedures and introduces the distributionally robust formulation based on BCOT. Section \ref{sec:convex-re} reformulates the minimax problem as a convex optimization problem. Section \ref{sec:example} presents two examples demonstrating the outperformance of our BCOT method. The code is publicly available at \url{https://github.com/hanbingyan/GaussianCOT_release}. Proofs are in the Appendix \ref{sec:append}.

We summarize the notations used in this paper. Let $I_n$ denote the $n \times n$ identity matrix. The trace of a matrix $X$ is denoted by $\text{tr}[X]$. The space of all symmetric matrices in $\R^{n \times n}$ is denoted by $\bS^n$. The symbols $\bS^n_+$ and $\bS^n_{++}$ represent the cones of symmetric positive semi-definite and positive definite matrices, respectively. For $X, Y \in \bS^n$, the relations $X \succeq Y$ and $X \succ Y$ indicate that $X - Y \in \bS^n_+$ and $X - Y \in \bS^n_{++}$, respectively. For $X \in \bS^n_+$, $\sqrt{X}$ denotes the square root matrix of $X$.

\section{Problem formulation}\label{sec:form}
\subsection{Kalman filter}
Consider a discrete-time linear stochastic state-space model with $t \in \mathbb{N}$:
\begin{align}
	x_t &= A_t x_{t-1} + w_t, \quad w_t \sim N(0, Q_t), \label{state} \\
	y_t &= C_t x_t + v_t, \quad v_t \sim N(0, R_t). \label{measure}
\end{align}
Here, $x_t$ is the $n$-dimensional unobserved state vector. The state matrix $A_t \in \R^{n \times n}$. The process noise vector $w_t \in \R^n$ follows a zero-mean Gaussian distribution with covariance $Q_t \in \bS^n_{++}$. Let $y_t$ be the $m$-dimensional observed measurement variable. The observation matrix $C_t \in \R^{m \times n}$. The measurement noise vector $v_t \in \R^m$ follows a zero-mean Gaussian distribution with covariance $R_t \in \bS^m_{++}$. We assume that the initial state $x_0 \sim N(\hat{x}_0, \Sigma_0)$. The noise processes $\{w_t\}$ and $\{v_t\}$ are white, uncorrelated, and independent of the initial state $x_0$.

The Kalman filter is a widely-used method for estimating the unobserved state $x_t$ based on noisy observations $y_{1:t}$. It operates recursively through two steps.  Suppose the conditional distribution of the state $x_{t-1}$ given measurements $y_{1:t-1}$ has been obtained as $\p(x_{t-1}| y_{1:t-1}) = N(\hat{x}_{t-1}, \Sigma_{t-1})$ with $\Sigma_{t-1} \in \bS^n_{++}$. 
\begin{enumerate}
	\item \textbf{Prediction Step}: Calculate $\p(x_t, y_t| y_{1:t-1})$ using the previous $\p(x_{t-1}| y_{1:t-1})$ and the transition $\p(x_t, y_t|x_{t-1})$ implied by the state-space model \eqref{state}-\eqref{measure}. Note that $\p(x_t, y_t | x_{t-1})$ is
	\begin{equation*}
		N\left(\begin{pmatrix}
			A_t x_{t-1} \\
			C_t A_t x_{t-1}
		\end{pmatrix}, \begin{pmatrix}
			Q_t, & Q_t C^\top_t \\
			C_t Q_t, & C_t Q_t C^\top_t + R_t
		\end{pmatrix}\right).
	\end{equation*}
	Therefore, $\p(x_t, y_t| y_{1:t-1})$ is still Gaussian. We obtain $\p(x_t, y_t| y_{1:t-1}) = N(\mu_t, V_t)$, where $\mu^\top_{t} = (\mu^\top_{t, x}, \mu^\top_{t, y})$ and
	\begin{equation*}
		V_{t} = \begin{pmatrix}
			V_{t, xx}, & V_{t, xy} \\
			V_{t, yx}, & V_{t, yy}
		\end{pmatrix},
	\end{equation*}
	with 
	\begin{equation}\label{Cov} 
		\begin{aligned}
			\mu_{t, x} = & A_t \hat{x}_{t-1}, \quad \mu_{t, y} = C_t A_t \hat{x}_{t-1}, \\
			V_{t, xx} = & 	A_t \Sigma_{t-1} A^\top_t + Q_t, \\
			V_{t, xy} = &   A_t \Sigma_{t-1} A^\top_t C^\top_t + Q_t C^\top_t, \\
			V_{t, yx} = &  C_t A_t \Sigma_{t-1} A^\top_t + C_t Q_t, \\
			V_{t, yy} = & C_t A_t \Sigma_{t-1} A^\top_t C^\top_t + C_t Q_t C^\top_t + R_t.
		\end{aligned}
	\end{equation}
	We assume $V_t \succ 0$ in this paper.
	
	\item \textbf{Update Step}: Calculate the posterior distribution with the new measurement $y_t$. This is done by considering the minimum mean square error (MMSE) problem:
	\begin{equation}\label{mmse}
		\inf_{f \in \cL} \E_\p \left[ \| x_t - f(y_t) \|^2 \Big| y_{1:t-1} \right], 
	\end{equation} 
	where $\cL$ denotes the set of all measurable functions; see \cite[eq.14]{zorzi2016robust} and \cite[eq.11]{nips2018kalman}. Note that the estimator depends on $y_{1:t}$, not just on $y_t$, but this dependency is simplified here. The expectation is under $\p(x_t, y_t| y_{1:t-1}) = N(\mu_t, V_{t})$ obtained in the prediction step. The MMSE estimator is the conditional expectation $ \E_\p[x_t | y_{1:t}]$. Furthermore, the Gaussian assumption yields $\p(x_t | y_{1:t}) = N(\hat{x}_{t}, \Sigma_{t})$ with 
	\begin{equation}\label{update} 
		\begin{aligned}
			\hat{x}_{t} =& V_{t, xy} V^{-1}_{t, yy} (y_t - \mu_{t, y}) + \mu_{t, x}, \\ \quad \Sigma_{t} =& V_{t, xx} - V_{t, xy} V^{-1}_{t, yy} V_{t, yx}.
		\end{aligned}
	\end{equation}
\end{enumerate}

\subsection{Distributionally robust framework}
Under the linear state-space model \eqref{state}-\eqref{measure} with known noise covariance matrices $Q_t$ and $R_t$, the Kalman filter is the optimal estimator in terms of the MMSE \eqref{mmse} when $\{w_t\}$ and $\{v_t\}$ are Gaussian, zero-mean, uncorrelated, and white. If $\{w_t\}$ and $\{v_t\}$ are not Gaussian, please refer to \cite[p.130]{simon2006optimal} and \cite[Theorem 2.1, p.47]{anderson1979} for a discussion on the properties of the Kalman filter. Additionally, \cite{idan2010cauchy,idan2013multivariate} demonstrate the optimality of the Cauchy estimator in a Cauchy environment. However, in practical scenarios, the user-specified  $Q_t$ and $R_t$ may not accurately represent the true covariance matrices, potentially resulting in model misspecification or even divergence in the Kalman filter. Additionally, the covariance $\Sigma_{t-1}$ of the previous estimate can also be misspecified and contribute to the estimation errors.

To enhance robustness, a commonly adopted approach involves considering alternative distributions similar to, but distinct from, the reference model. In this study, we explore alternative models that remain linear and Gaussian, as expressed by
\begin{align}
	\bar{x}_t &= A_t \bar{x}_{t-1} + \bar{w}_t, \quad \bar{w}_t \sim N(0, \bar{Q}_t), \label{alt-state} \\
	\bar{y}_t &= C_t \bar{x}_t + \bar{v}_t, \quad \bar{v}_t \sim N(0, \bar{R}_t). \label{alt-measure}
\end{align}
Unlike the approach in \cite{nips2018kalman}, we assume no uncertainty exists in the state matrix $A_t$ and observation matrix $C_t$. This assumption is valid in some applications, such as target tracking and pairs trading, where $A_t$ and $C_t$ are easy to obtain.

The alternative model \eqref{alt-state}-\eqref{alt-measure} suggests the transition kernel $\barp(\bar{x}_t, \bar{y}_t|\bar{x}_{t-1})$ as follows:
\begin{equation*}
	N\left(\begin{pmatrix}
		A_t \bar{x}_{t-1} \\
		C_t A_t \bar{x}_{t-1}
	\end{pmatrix},\begin{pmatrix}
		\bar{Q}_t, & \bar{Q}_t C^\top_t \\
		C_t \bar{Q}_t, & C_t \bar{Q}_t C^\top_t + \bar{R}_t
	\end{pmatrix}\right).
\end{equation*}
To maintain tractability, we assume that the alternative estimate employs the same mean $\hat{x}_{t-1}$, while allowing for a potentially different covariance estimate. Hence, we define $\barp(\bar{x}_{t-1}| y_{1:t-1}) = N(\hat{x}_{t-1}, \bar{\Sigma}_{t-1})$. It is important to note that the uncertainty in $\hat{x}_{t-1}$ is handled in the update step at time $t-1$. Moving to the current time $t$, the framework should aim to determine a robust $\hat{x}_t$, thus maintaining the recursive nature of the problem.

At a given time $t$, our approach considers uncertainty not only in the noise covariance matrices $Q_t$ and $R_t$, but also in the previous $\Sigma_{t-1}$. This enables us to correct estimation errors from the previous step, preventing their propagation and accumulation in future steps. This dynamic approach, spanning two time steps, differs from the static, single-time-step framework in the work by \cite{nips2018kalman}. Therefore, our optimization problem \eqref{eq:minimax} involves finding robust values for $\bar{Q}_t$, $\bar{R}_t$, and $\bar{\Sigma}_{t-1}$ used in the update step, instead of just $\bar{Q}_t$ and $\bar{R}_t$.

A technical challenge is quantifying the distance between the nominal model $ \p_{t-1} := \p(x_{t-1}, x_t, $ $y_t|y_{1:t-1})$ and an alternative model $\barp_{t-1} := \barp(\bar{x}_{t-1}, \bar{x}_t, \bar{y}_t | y_{1:t-1})$. Two methods lacking temporal structure have been proposed in the literature. \cite{zorzi2016robust} employed a KL divergence-based approach, while \cite{nips2018kalman} utilized the OT paradigm. KL divergence is asymmetric and requires absolute continuity for finite values, whereas the OT framework does not have these restrictions. We provide a brief overview of the OT framework below.

A coupling, denoted by $\pi$, is a Borel probability measure that has $\p_{t-1}$ and $\barp_{t-1}$ as its marginals. Denote $\Pi(\p_{t-1}, \barp_{t-1})$ as the set of all such couplings. $\pi(x_{t-1:t}, y_t, \bar{x}_{t-1:t}, \bar{y}_t)$ is commonly termed a transport plan between $\p_{t-1}$ and $\barp_{t-1}$. In simpler terms, one can interpret $\pi(x_{t-1:t}, y_t, \bar{x}_{t-1:t}, \bar{y}_t)$ as the probability mass moved from $(x_{t-1:t}, y_t)$ to $(\bar{x}_{t-1:t}, \bar{y}_t)$.

We suppose transporting one unit of mass from $(x_{t-1:t}, y_t)$ to $(\bar{x}_{t-1:t}, \bar{y}_t)$ incurs a cost of $c(x_{t-1:t}, \bar{x}_{t-1:t}, y_t, \bar{y}_t)$. The classical OT problem is formulated as 
\begin{equation}\label{eq:OT-dist}
	\begin{aligned}
		& \cW(\p_{t-1}, \barp_{t-1}) := \inf_{ \pi \in \Pi(\p_{t-1}, \barp_{t-1})} \int c(x_{t-1:t}, \bar{x}_{t-1:t}, y_t, \bar{y}_t) d\pi.
	\end{aligned}
\end{equation}
Here, $\cW(\p_{t-1}, \barp_{t-1})$ quantifies the difficulty of reshaping the measure $\p_{t-1}$ to $\barp_{t-1}$. After observing $y_{t-1}$, which is also $\bar{y}_{t-1}$, we assume a quadratic cost:
\begin{align*}
	c(x_{t-1:t}, \bar{x}_{t-1:t}, y_t, \bar{y}_t) = & \|x_{t-1} - \bar{x}_{t-1}\|^2 + \| x_t - \bar{x}_t \|^2 + \| y_t - \bar{y}_t \|^2,
\end{align*}
where $y_t$ and $\bar{y}_t$ are undisclosed and random at time $t-1$. To calculate $\cW(\p_{t-1}, \barp_{t-1})$, we establish the following lemma, which extends \cite{givens1984class} slightly.
\begin{lemma}\label{cor:Gaussian}
	Suppose two $n$-dimensional multivariate normal random vectors are given by $\xi_1 \sim \p_1 = N(m_1, M_1)$ and $\xi_2 \sim \p_2 = N(m_2, M_2)$, where $M_1, M_2 \in \bS^n_+$. Denote $P, S, U$ as symmetric matrices with suitable dimensions. Then
	\begin{align*}
		& \inf_{ \pi \in \Pi(\p_1, \p_2)} \int (\xi_1^\top P \xi_1 + \xi_2^\top S \xi_2 + \xi_1^\top U \xi_2) d\pi \\
		& \qquad = m^\top_1 P m_1 + m^\top_2 S m_2 + m^\top_1 U m_2 \\
		&\qquad \quad + \tr[PM_1] + \tr[SM_2] - \tr \left[\sqrt{\sqrt{M_1} U M_2 U \sqrt{M_1}} \right].
	\end{align*} 
\end{lemma}

To apply Lemma \ref{cor:Gaussian}, we note that $\p_{t-1}$ and $\barp_{t-1}$ are normal distributions. Specifically, $\p(x_{t-1}, x_t, y_t| y_{1:t-1})$ is a normal distribution with the mean $(\hat{x}^\top_{t-1}, (A_t \hat{x}_{t-1})^\top, (C_t A_t \hat{x}_{t-1})^\top)^\top$ and covariance
\begin{equation}\label{eq:Pt-1}
	\begin{pmatrix}
		\Sigma_{t-1}, & \Sigma_{t-1} A^\top_t, & \Sigma_{t-1} A^\top_t C^\top_t \\ 
		A_t \Sigma_{t-1}, & J_{22}, & J_{23} \\
		C_t A_t \Sigma_{t-1}, & J^\top_{23}, & J_{33}
	\end{pmatrix},
\end{equation}
where
\begin{align*}
	J_{22} =& A_t \Sigma_{t-1} A^\top_t + Q_t, \; J_{23} = A_t \Sigma_{t-1} A^\top_t C^\top_t + Q_t C^\top_t, \\
	J_{33} = &  C_t A_t \Sigma_{t-1} A^\top_t C^\top_t + C_t Q_t C^\top_t + R_t.
\end{align*}
The distribution $\barp_{t-1}$ is obtained similarly.

Next, the OT framework inspired by \cite{nips2018kalman}, but with fixed $A_t$ and $C_t$, examines alternative models within a Wasserstein ball of radius $\varepsilon > 0$:
\begin{align*}
	\cV_{\varepsilon, t - 1} := \big\{ \barp_{t-1} : & \cW (\p_{t-1},\barp_{t-1}) \leq \varepsilon, \\
	& \barp_{t-1} \text{ is defined by \eqref{alt-state}-\eqref{alt-measure}} \big\}.
\end{align*}

The distributionally robust MMSE problem in the OT framework with fixed $A_t$ and $C_t$ is expressed as follows:
\begin{equation}\label{eq:OTminimax}
	\begin{aligned}
		& \inf_{f \in \cL} \sup_{\barp_{t-1}} \E_{\barp_{t-1}} \left[ \| \bar{x}_t - f(\bar{y}_t) \|^2 \Big| y_{1:t-1} \right] \\
		& \textrm{ subject to} \quad \barp_{t-1} \in \cV_{\varepsilon, t - 1}, \; \bar{Q}_t \in \bS^n_+, \, \bar{\Sigma}_{t-1} \in \bS^n_+, \\
		& \qquad \qquad \text{ and } \bar{R}_t \succeq \delta I_m,
	\end{aligned}
\end{equation}
where $\delta > 0$ is fixed and sufficiently small. The constraint $\bar{R}_t \succeq \delta I_m$ ensures that the problem does not degenerate. The goal of \eqref{eq:OTminimax} is to determine robust values for $\bar{Q}_t$, $\bar{R}_t$, and $\bar{\Sigma}_{t-1}$ to be used in the update step.

\subsection{Bicausal optimal transport}
To justify our BCOT formulation, we highlight two shortcomings of the classical OT formulation \eqref{eq:OTminimax}. First, the covariance matrices in $\p_{t-1}$ and $\barp_{t-1}$ have dimensions $(2n + m, 2n + m)$, resulting in increased computational complexity. Second, since our problem spans two time steps, a fundamental requirement for the transport plan $\pi$ on temporal data is the non-anticipative condition \cite{lassalle2013causal}. We present this condition when $y_{1:t-1}$ is observed, assuming $\bar{y}_{1:t-1} = y_{1:t-1}$. Given $(x_{t-1},  y_{1:t-1})$, $(\bar{x}_{t-1}, \bar{y}_{t-1})$ should be independent of $(x_t,  y_t)$ under $\pi$. Mathematically, a transport plan $\pi$ should satisfy
\begin{equation}\label{eq:causal}
	\begin{aligned}
		& \pi (d\bar{x}_{t-1}, d\bar{y}_{t-1} | x_{t-1}, y_{1:t-1}, x_t, y_t) =\pi (d\bar{x}_{t-1}, d\bar{y}_{t-1} | x_{t-1}, y_{1:t-1}), \; \pi\text{-a.s.}
	\end{aligned}
\end{equation}
This condition \eqref{eq:causal} is known as the causality condition, and a transport plan satisfying \eqref{eq:causal} is called {\it causal} by \cite{lassalle2013causal}. If the same condition also holds when the positions of $(x, y)$ and $(\bar{x}, \bar{y})$ are exchanged, i.e.
\begin{equation}\label{eq:anticausal}
	\begin{aligned}
		& \pi (dx_{t-1}, dy_{t-1} | \bar{x}_{t-1}, \bar{y}_{1:t-1}, \bar{x}_t, \bar{y}_t) = \pi (dx_{t-1}, dy_{t-1} | \bar{x}_{t-1}, \bar{y}_{1:t-1}), \; \pi\text{-a.s.},
	\end{aligned}
\end{equation}
the transport plan is called {\it bicausal}. Denote $\Pi_{bc}(\p_{t-1}, \barp_{t-1})$ as the set of all bicausal transport plans between $\p_{t-1}$ and $\barp_{t-1}$. \eqref{eq:causal}-\eqref{eq:anticausal} only consider the time step $t-1$ since our formulation spans two steps, with the case at time $t$ automatically satisfied. For consideration of more than two steps, the general version of \eqref{eq:causal}-\eqref{eq:anticausal} should be employed, as in \cite{lassalle2013causal}.

The bicausal or adapted Wasserstein distance is defined as
\begin{align}\label{eq:bicau-dist}
	& \cW_{bc} (\p_{t-1}, \barp_{t-1})   := \inf_{ \pi \in \Pi_{bc}(\p_{t-1}, \barp_{t-1})} \int c(x_{t-1:t}, \bar{x}_{t-1:t}, y_t, \bar{y}_t) d \pi,
\end{align}
where the integral is over $x_{t-1:t}, \bar{x}_{t-1:t}, y_t, \bar{y}_t$.

Using Lemma \ref{cor:Gaussian}, Lemma \ref{cor:Wasser} gives an explicit solution to the bicausal Wasserstein distance \eqref{eq:bicau-dist}. The highest dimension of the matrices is reduced to $(n + m, n + m)$ in the bicausal Wasserstein distance formula. Notably, \eqref{eq:bicau-dist} can be computed recursively using dynamic programming principle; see \cite[Proposition 5.2]{backhoff2017causal}.  
\begin{lemma}\label{cor:Wasser}
	Denote the following matrices
	\begin{equation*}
		\sigma^2(Q_t, R_t) := \begin{pmatrix}
			Q_t & Q_t C^\top_t \\
			C_t Q_t & C_t Q_t C^\top_t + R_t
		\end{pmatrix},
	\end{equation*}
	\begin{equation*}
		\sigma^2(\bar{Q}_t, \bar{R}_t) := \begin{pmatrix}
			\bar{Q}_t & \bar{Q}_t C^\top_t \\
			C_t \bar{Q}_t & C_t \bar{Q}_t C^\top_t + \bar{R}_t
		\end{pmatrix},
	\end{equation*}
	and $H := I_n + A^\top_t A_t + A^\top_t C^\top_t C_t A_t$. Then the bicausal Wasserstein distance \eqref{eq:bicau-dist} is a function of $\bar{Q}_t, \bar{R}_t, \bar{\Sigma}_{t-1}$ given by $\cW_{bc} (\p_{t-1}, \barp_{t-1}) = w(\bar{Q}_t, \bar{R}_t, \bar{\Sigma}_{t-1})$, where
	\begin{equation}\label{w}
		\begin{aligned}
			& w(\bar{Q}_t, \bar{R}_t, \bar{\Sigma}_{t-1}) \\
			&:= \tr \Big[\sigma^2(Q_t, R_t) + \sigma^2(\bar{Q}_t, \bar{R}_t)  \\ 
			&\qquad \quad - 2 \sqrt{\sqrt{\sigma^2(Q_t, R_t)} \sigma^2(\bar{Q}_t, \bar{R}_t) \sqrt{\sigma^2(Q_t, R_t)} } \Big] \\ 
			& \quad  + \tr \Big[ H \Sigma_{t-1} + H \bar{\Sigma}_{t-1} - 2 \sqrt{\sqrt{\Sigma_{t-1}} H \bar{\Sigma}_{t-1} H \sqrt{\Sigma_{t-1}} } \Big].
		\end{aligned}
	\end{equation}
\end{lemma}

Denote the bicausal Wasserstein ball with a fixed radius $\varepsilon$ as
\begin{align*}
	\cB_{\varepsilon, t - 1} := \big\{ \barp_{t-1} : & \cW_{bc} (\p_{t-1},\barp_{t-1}) \leq \varepsilon,\; \barp_{t-1} \text{ is defined by \eqref{alt-state}-\eqref{alt-measure}} \big\}.
\end{align*}
The effect of the radius $\varepsilon$ is highlighted by the numerical analysis in Section \ref{sec:example}.

In contrast to the OT formulation \eqref{eq:OTminimax}, the BCOT approach considers the bicausal Wasserstein ball instead:
\begin{equation}\label{eq:minimax}
	\begin{aligned}
		& \inf_{f \in \cL} \sup_{\barp_{t-1}} \E_{\barp_{t-1}} \left[ \| \bar{x}_t - f(\bar{y}_t) \|^2 \Big| y_{1:t-1} \right] \\
		& \textrm{ subject to} \quad \barp_{t-1} \in \cB_{\varepsilon, t - 1}, \; \bar{Q}_t \in \bS^n_+, \, \bar{\Sigma}_{t-1} \in \bS^n_+, \\
		& \qquad \qquad \text{ and } \bar{R}_t \succeq \delta I_m,
	\end{aligned}
\end{equation}
where $\delta > 0$ is fixed and sufficiently small. The main distinction between the OT and BCOT methods, as shown in equations \eqref{eq:OTminimax} and \eqref{eq:minimax} respectively, lies in the inclusion of the bicausality constraint \eqref{eq:causal}-\eqref{eq:anticausal}. Additionally, it is worth mentioning that \cite{yang2022decision} explored a relevant problem with COT and general distributions, instead of the BCOT and Gaussian setting in this work.

In the next section, we reduce \eqref{eq:minimax} to a convex problem.

\section{Convex reformulation}\label{sec:convex-re}
\subsection{Minimax theorem}
While the formulation \eqref{eq:minimax} is non-convex and infinite dimensional in general, Theorem \ref{thm:covex} proves that \eqref{eq:minimax} is equivalent to a nonlinear semi-definite programming problem.
\begin{theorem}\label{thm:covex}
	Denote the following matrices
	\begin{align*}
		K & := C_t A_t \bar{\Sigma}_{t-1} A^\top_t C^\top_t + C_t \bar{Q}_t C^\top_t + \bar{R}_t, \\
		M & := C_t A_t \bar{\Sigma}_{t-1} A^\top_t + C_t \bar{Q}_t.
	\end{align*}
	Define
	$$ F(\bar{Q}_t, \bar{R}_t, \bar{\Sigma}_{t-1})  := \tr[ A_t \bar{\Sigma}_{t-1} A^\top_t + \bar{Q}_t] - \tr[ M M^\top K^{-1}]. 
	$$
	Then the distributionally robust MMSE problem \eqref{eq:minimax} is equivalent to the following optimization problem
	\begin{equation}\label{eq:main-pro} 
		\begin{aligned}
			& \sup_{\bar{Q}_t, \bar{R}_t, \bar{\Sigma}_{t-1}}  F(\bar{Q}_t, \bar{R}_t, \bar{\Sigma}_{t-1}) \\
			& \textrm{s.t.} \; w(\bar{Q}_t, \bar{R}_t, \bar{\Sigma}_{t-1}) \leq \varepsilon, \\
			& \quad \bar{Q}_t \succeq 0, \; \bar{R}_t \succeq \delta I_m, \; \bar{\Sigma}_{t-1} \succeq 0,
		\end{aligned}
	\end{equation}
	with the function $w(\cdot)$ defined in \eqref{w}. Moreover, if $(Q_*, R_*, \Sigma_*)$ is an optimizer for \eqref{eq:main-pro}, then the robust estimator in \eqref{eq:minimax} is given by
	\begin{equation}\label{rob_f}
		f(\bar{y}_t) = M^\top_* K^{-1}_* \bar{y}_t + A_t \hat{x}_{t-1} - M^\top_* K^{-1}_* C_t A_t \hat{x}_{t-1},
	\end{equation}
	where $M_*$ and $K_*$ are defined with $(Q_*, R_*, \Sigma_*)$. The worst-case distribution $\barp_*$ is determined by replacing the corresponding terms in \eqref{update} and \eqref{Cov} with $(Q_*, R_*, \Sigma_*)$.
\end{theorem}

In general, the optimizer $\Sigma_* \neq \Sigma_{t-1}$, confirming our assertion that the previous covariance matrix is adjusted. One might wonder if $\cW_{bc} (\p_{t-1}, \barp_{t-1})$ could be replaced by the classical OT metric \cite{nips2018kalman} or the KL divergence \cite{zorzi2016robust}. However, these methods lack the bicausality constraint \eqref{eq:causal}-\eqref{eq:anticausal}. Consequently, they might determine $\bar{\Sigma}_{t-1}$ in an anticipative way. When the problem spans two time steps $t-1$ and $t$, it is essential to recognize the evolution of information flows by enforcing the bicausality condition.

One advantage of our framework is the optimization problem \eqref{eq:main-pro} is indeed convex, as proved in Lemma \ref{lem:convex}. The proof is lengthy but straightforward.
\begin{lemma}\label{lem:convex}
	The optimization problem \eqref{eq:main-pro} is convex.
\end{lemma}

\subsection{Optimization with $LDL^\top$ decomposition}
The optimization problem \eqref{eq:main-pro} involves a nonlinear objective function and nonlinear semi-definite inequality constraints. There exists a substantial body of literature addressing nonlinear semi-definite programming, with specialized algorithms developed for this purpose. This paper adopts a general methodology outlined in \cite{benson2003solving} to reformulate the semi-definite constraints as convex constraints.

\cite[Theorem 1]{benson2003solving} establishes that the interior of $\bS^n_+$ is $\bS^n_{++}$. For any $\bar{Q} \in \bS^n_+$, there exists a unit triangular matrix $L$ and a unique diagonal matrix $D$ such that $\bar{Q} = L D L^T$, a decomposition known as the $L D L^T$ decomposition. \cite{benson2003solving} provides a useful method for computing the entries of $D$. For $j \in \{ 1, ... , n\}$, we denote the blocks of $\bar{Q}$ as
\begin{equation*}
	\bar{Q} = 
	\begin{pNiceMatrix}
		Z & c & * \\
		c^\top & b & * \\
		* & * & *
	\end{pNiceMatrix},
\end{equation*} 
where $Z$ is the $(j-1) \times (j-1)$ principle submatrix, $c$ is a column vector with $j-1$ elements, and $b$ is the $(j, j)$ entry. \cite[Theorem 2]{benson2003solving} proves that if $\bar{Q} \in \bS^n_{++}$, then the matrix $Z$ is nonsingular, and the $j$-th diagonal element of $D$ is given by $d_j (\bar{Q}) = b - c^\top Z^{-1} c$. \cite[Theorem 3]{benson2003solving} shows that $d_j (\bar{Q})$ is concave. Consequently, we can replace the positive semi-definite constraint $\bar{Q} \succeq 0$ with $d_j (\bar{Q}) \geq 0, j = 1, ... , n$. Similarly, positive semi-definite constraints on $\bar{\Sigma}$ can be replaced. To ensure positive definiteness of $\bar{R}$, we can impose $d_j (\bar{R}) \geq \underline{\delta}, \, j = 1, ... , m$, where $\underline{\delta} > 0$ is small. In fact, for numerical stability during implementation, positive semi-definite constraints on $\bar{Q}$ and $\bar{\Sigma}$ are also replaced with small positive lower bounds on all $d_j$, as detailed in \cite[Section 6]{benson2003solving}.   

Using the $LDL^\top$ decomposition, we transform the semi-definite constraints into convex constraints, allowing us to solve the problem with an interior-point method. Our implementation uses the {\it minimize} function from the {\it scipy.optimize} package in Python. To ensure that intermediate solutions remain feasible, we set the argument $``keep\_feasible=True"$ in the $NonlinearConstraint$ function. Although convergence may require many iterations, our experiments show that early stopping after a few iterations can still yield good performance. For computational efficiency, we limit the maximum number of iterations to 20 in all OT and BCOT experiments.

Algorithm \ref{algo:BCOT} summarizes the distributionally robust Kalman filter utilizing the BCOT method as follows:
\begin{algorithm}[h]
	\small
	\begin{algorithmic}[1]
		\State {\bf Input:} State matrix $A_t \in \R^{n \times n}$; Observation matrix $C_t \in \R^{m \times n}$; Covariance matrix $\Sigma_{t-1} \in \bS^n_{++}$; State estimate $\hat{x}_{t-1}$; BCOT radius $\varepsilon > 0$; Threshold $\delta > 0$ in \eqref{eq:main-pro}
		\State {\bf Prediction:} Construct the nominal distribution $\p_{t-1} = \p(x_{t-1}, x_t, y_t| y_{1:t-1})$ using \eqref{eq:Pt-1}
		\State {\bf Observation:} Observe $\bar{y}_t$, which is also $y_t$ 
		\State {\bf Update:} Use the interior-point method with the $LDL^\top$ decomposition to solve \eqref{eq:main-pro} and find an optimizer $(Q_*, R_*, \Sigma_*)$
		\State Define $V^*_t$ as in \eqref{Cov}, with $(Q_t, R_t, \Sigma_{t-1})$ replaced by $(Q_*, R_*, \Sigma_*)$
		\State {\bf Output:} $\hat{x}_{t} = V^*_{t, xy} (V^*_{t, yy})^{-1} (\bar{y}_t - C_t A_t \hat{x}_{t-1}) + A_t\hat{x}_{t-1} = f(\bar{y}_t)$ in \eqref{rob_f}
		\State $\Sigma_{t} = V^*_{t, xx} - V^*_{t, xy} (V^*_{t, yy})^{-1} V^*_{t, yx}$
	\end{algorithmic}
	\caption{Distributionally robust Kalman filter at time $t$}
	\label{algo:BCOT}
\end{algorithm}

\section{Examples}\label{sec:example}
\subsection{Target tracking}
We examine the target tracking simulation discussed in \cite{huang2017novel}. Suppose the target moves in a two-dimensional space under continuous white noise acceleration. We observe noisy measurements of the target's position, denoted as $y_t \in \R^2$. The state vector is $x_t = [x_{1, t}, x_{2, t}, \dot{x}_{1, t}, \dot{x}_{2, t}]^\top$, where $(x_{1, t}, x_{2, t})$ represents the Cartesian coordinate of the target, and $(\dot{x}_{1, t}, \dot{x}_{2, t})$ is the corresponding velocity vector. The state matrix $A_t$ and the observation matrix $C_t$ are as follows:
\begin{equation*}
	A_t = \begin{pNiceMatrix}[r]
		I_2 & \Delta t I_2  \\
		0 & I_2
	\end{pNiceMatrix}, \quad C_t = \begin{pNiceMatrix}[l]
		I_2 & 0
	\end{pNiceMatrix},
\end{equation*}
where $\Delta t = 1s$ denotes the sampling interval. The noise distribution setup mirrors that of \cite{huang2017novel}. The real state noise covariance matrix $Q^0_t$ and the observation noise covariance matrix $R^0_t$ are given by
\begin{equation*}
	Q^0_t = q \left[6.5 + 0.5 \cos \left(\frac{\pi t}{T} \right) \right] \begin{pNiceMatrix}[r]
		\frac{(\Delta t)^3}{3} I_2 & \frac{(\Delta t)^2}{2} I_2  \\
		\frac{(\Delta t)^2}{2} I_2 & \Delta t I_2
	\end{pNiceMatrix}
\end{equation*}
and
\begin{equation*}
	R^0_t = r \left[0.1 + 0.05 \cos\left(\frac{\pi t}{T} \right) \right] \begin{pNiceMatrix}[r]
		1 & 0.5  \\
		0.5 & 1
	\end{pNiceMatrix},
\end{equation*}
respectively. Here, $T = 100s$ represents the simulation time. We set $q = 10 m^2/s^3$ and $r = 50 m^2$. 

In the simulation, the agent needs to specify the noise covariance matrices. The procedure is as follows: Assuming the agent observes training samples within a significantly shorter simulation duration $T' = 10$, the agent employs the EM algorithm to estimate static $Q$ and $R$ matrices, which are time-independent. Evidently, the estimated constant covariance matrices are misspecified.

We assess the effectiveness of four filtering methods:
\begin{enumerate}
	\item EM: It employs the EM algorithm with constant covariance matrices calibrated from training samples. It serves as an out-of-sample test for model misspecification and does not incorporate robustness. Additionally, our framework can accommodate other Kalman-like filters as reference models.
	\item KL: The second method is based on the KL divergence, extending the approach detailed in \cite{zorzi2016robust} to find robust values of $\bar{Q}_t$, $\bar{R}_t$, and $\bar{\Sigma}_{t-1}$. 
	\item OT: The third method adopts the OT formulation in \eqref{eq:OTminimax}. 
	\item BCOT: The fourth method is our bicausal formulation \eqref{eq:minimax}.
\end{enumerate} 
In all methods, the matrices $A_t$ and $C_t$ are known. The radius of the Wasserstein ball significantly influences performance. We evaluate radii ranging from $0.5$ to $4.0$ in increments of $0.5$. Each radius undergoes simulation across $I = 10$ instances. To evaluate the precision of the filtered outcome, we introduce the Root Mean Square Error (RMSE) for position measurements at time $t$ of instance $i$ as
\begin{equation*}
	\text{RMSE}_{i, t} = \sqrt{\Big[ (x^i_{1, t} - \hat{x}^i_{1, t})^2 + (x^i_{2, t} - \hat{x}^i_{2, t})^2 \Big]}.
\end{equation*}
Here, $(\hat{x}^i_{1, t}, \hat{x}^i_{2, t})$ denotes the estimated coordinate in the $i$th simulation at time $t$, while $(x^i_{1, t}, x^i_{2, t})$ represents the actual coordinate. The RMSE for velocity measurements is defined similarly. 

\begin{table*}
	\footnotesize
	\centering
	\begin{tabular}{ccccccccc}
		\hline
		Radius & 0.5 & 1.0 & 1.5 & 2.0 & 2.5 & 3.0 & 3.5 & 4.0 \\
		\hline 
		Mean (BCOT-EM) & $-0.0842$ & $-0.3048$ & $-0.6042$ & $-0.4328$ & $-0.66$ &   $-0.3986$ &  $-0.2088$ & $-0.434$ \\
		\hline 
		STD(BCOT-EM) & 0.5772 & 1.3232 & 1.8824 & 1.4739 & 2.0594 & 1.187 & 0.9878 & 1.5237 \\
		\hline 
		BCOT Mean/STD ratio & $-0.1459$ & $-0.2303$ & $-0.321$ &  $-0.2936$ & $-0.3205$ & $-0.3358$ & $-0.2114$ & $-0.2849$ \\
		\hline
		Mean(OT-EM) & $-0.0489$ & $-0.295$ & $-0.5621$ & $-0.3366$ & $-0.6304$ & $-0.3382$ & $-0.1012$ & $-0.4093$ \\
		\hline
		STD(OT-EM) & 0.7131 & 1.4282 & 2.0356 & 1.5213 & 2.3487 & 1.3486 & 1.1658 & 1.808 \\
		\hline 
		OT Mean/STD ratio & $-0.0686$ & $-0.2065$ & $-0.2761$ & $-0.2212$ & $-0.2684$ & $-0.2508$ & $-0.0868$ & $-0.2264$ \\
		\hline
		Mean(KL-EM) & $-0.0753$ & $-0.2099$ & $-0.2605$ & $-0.1563$ & $-0.2592$ & $-0.142$ &  $-0.1049$ & $-0.1832$ \\
		\hline 
		STD(KL-EM) & 0.529  & 1.0534 & 1.0161 & 0.6312 & 0.8226 & 0.4826 & 0.549  & 0.8111\\
		\hline 
		KL Mean/STD ratio & $-0.1424$ & $-0.1993$ & $-0.2564$ & $-0.2476$ & $-0.3151$ & $-0.2942$ & $-0.191$ &  $-0.2259$ \\
		\hline
	\end{tabular}
	\caption{Statistics of the difference in RMSE of position}\label{tab:poi}
\end{table*}

\begin{table*}
	\footnotesize
	\centering
	\begin{tabular}{ccccccccc}
		\hline
		Radius & 0.5 & 1.0 & 1.5 & 2.0 & 2.5 & 3.0 & 3.5 & 4.0 \\
		\hline 
		Mean(BCOT-EM) & $-0.4664$ & $-2.7545$ & $-0.9604$ & $-1.161$ &  $-0.9274$ & $-0.4745$ & $-1.727$ & $-4.571$ \\
		\hline 
		STD(BCOT-EM) & 2.186 & 5.2453 & 3.5421 & 3.0911 & 2.5457 & 2.372  & 4.3455 & 14.3121 \\
		\hline 
		BCOT Mean/STD ratio & $-0.2134$ & $-0.5251$ & $-0.2711$ & $-0.3756$ & $-0.3643$ & $-0.2001$ & $-0.3974$ & $-0.3194$ \\
		\hline 
		Mean(OT-EM) & $-0.3777$ & $-2.9552$ & $-0.8091$ & $-0.9999$ & $-0.8515$ & $-0.2806$ & $-1.657$ &  $-4.5453$ \\
		\hline 
		STD(OT-EM) & 2.4186 & 5.6876 & 3.7311 & 3.4357 & 2.8709 & 2.6988 &  4.835  & 14.8376 \\
		\hline 
		OT Mean/STD ratio & $-0.1562$ & $-0.5196$ & $-0.2168$ & $-0.291$ & $-0.2966$ & $-0.104$ &  $-0.3427$ & $-0.3063$ \\
		\hline
		Mean(KL-EM) & $-0.3132$ & $-1.5354$ & $-0.5032$ & $-0.4153$ & $-0.3616$ & $-0.2439$ & $-0.6004$ & $-1.6875$ \\
		\hline 
		STD(KL-EM) & 1.8829 & 3.9144 & 2.8328 & 1.8739 & 1.6236 & 1.6588 & 2.2738 & 5.6053 \\
		\hline 
		KL Mean/STD ratio & $-0.1663$ & $-0.3922$ & $-0.1776$ & $-0.2216$ & $-0.2227$ & $-0.1471$ & $-0.2641$ & $-0.301$ \\
		\hline
	\end{tabular}
	\caption{Statistics of the difference in RMSE of velocity}\label{tab:vel}
\end{table*}

Tables \ref{tab:poi} and \ref{tab:vel} report the RMSE and provide the mean and standard deviation (STD) of the RMSE differences across all instances and time steps. For example, ``BCOT-EM" denotes the difference of RMSEs between BCOT and EM methods. Compared to the non-robust EM benchmark, the KL, OT, and BCOT methods exhibit reduced RMSE values for both position and velocity measurements, as reflected by negative average differences. The Mean/STD ratio stands for the ratio between the mean and STD of RMSE differences, representing the balance between error reduction and stability preservation.  Notably, the BCOT method demonstrates the best Mean/STD ratio and generally yields smaller RMSE values across various radii. 
\begin{figure}
	\centering
	\begin{minipage}{0.4\textwidth}
		\centering
		\includegraphics[width=0.8\textwidth]{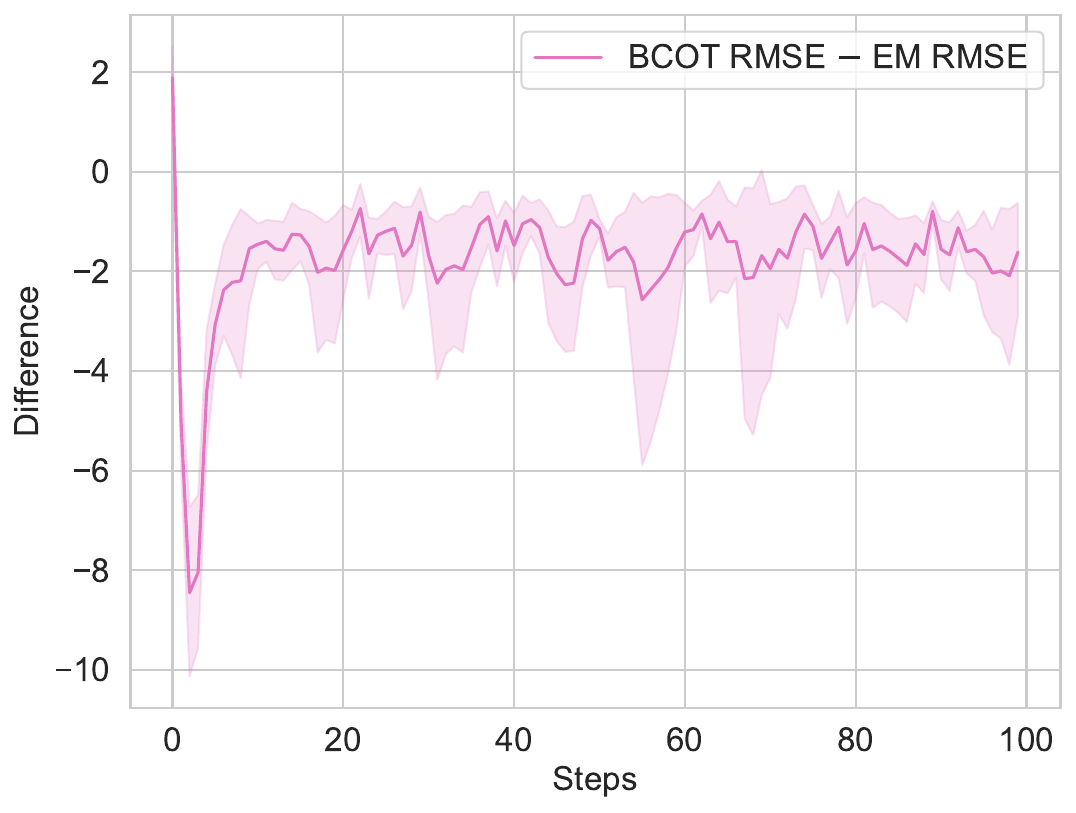}
		\subcaption{BCOT RMSE $-$ EM RMSE}\label{fig:EM}
	\end{minipage}
	\begin{minipage}{0.4\textwidth}
		\centering
		\includegraphics[width=0.8\textwidth]{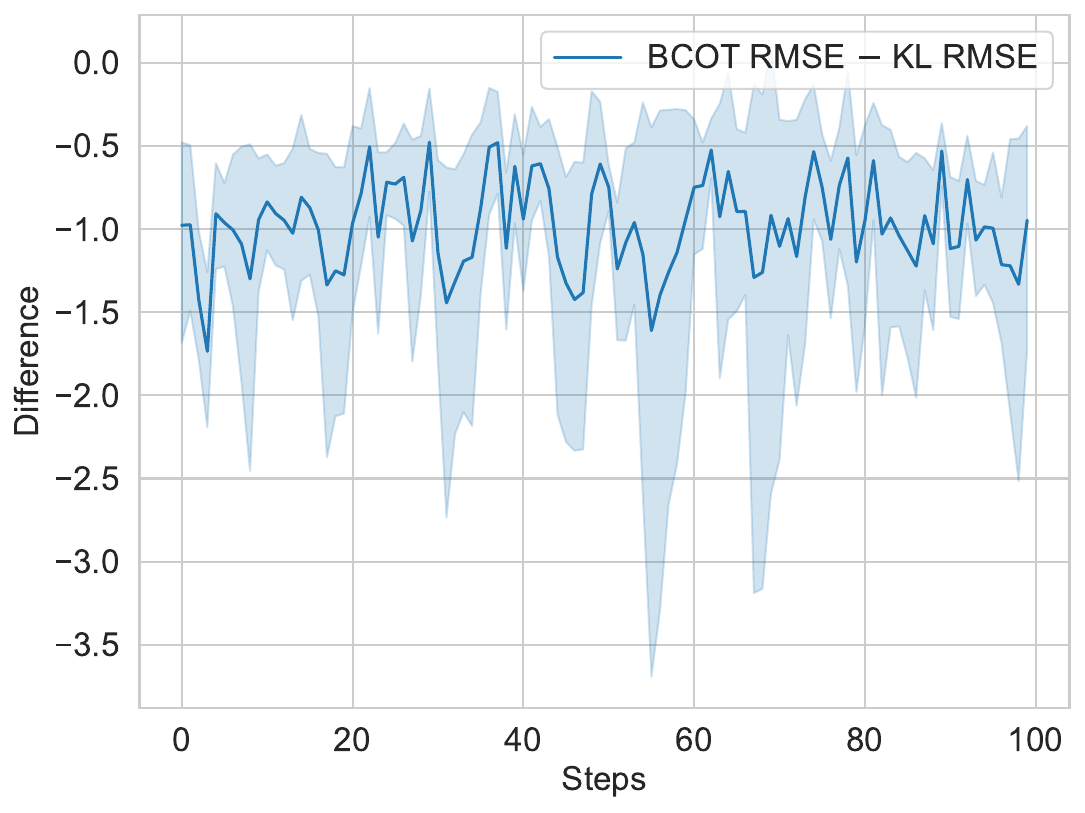}
		\subcaption{BCOT RMSE $-$ KL RMSE}\label{fig:KL}
	\end{minipage}
	\begin{minipage}{0.4\textwidth}
		\centering
		\includegraphics[width=0.7\textwidth]{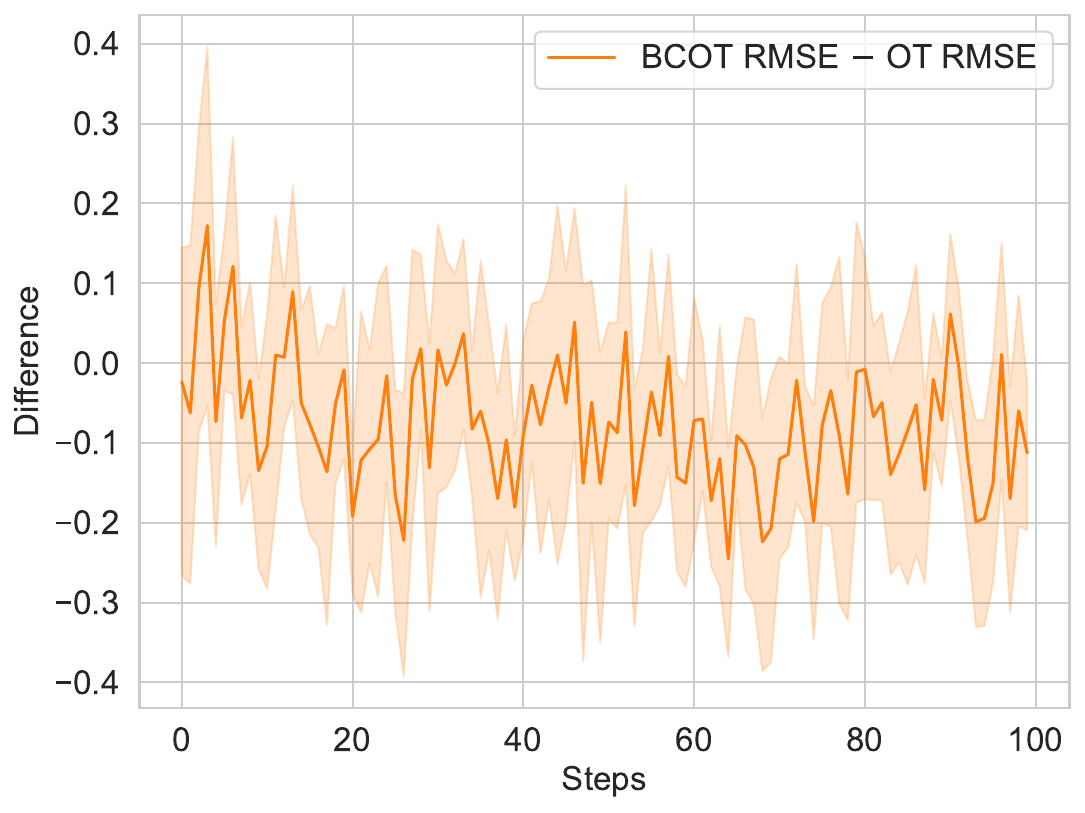}
		\subcaption{BCOT RMSE $-$ OT RMSE}\label{fig:OT}
	\end{minipage}%
	\caption{Time series of the difference in RMSE values}\label{fig:err}
\end{figure}

Figure \ref{fig:err} illustrates the RMSE differences among the four methods over time. The solid curve represents the mean value, while the shaded area indicates the 95\% confidence intervals, aggregated across all radii. As the confidence intervals consistently remain below zero, the BCOT method consistently outperforms the EM and KL methods. Furthermore, compared to OT, the BCOT method maintains a lower average RMSE, particularly in later time steps.

\subsection{Pairs trading}
In finance, a classical application of Kalman filters is pairs trading. This strategy involves analyzing two stock prices, $y_{1,t}$ and $y_{2,t}$, which exhibit high correlation due to shared characteristics. When one stock outperforms and the other underperforms relatively, the pairs trading entails buying the underperforming stock and selling the outperforming one, with the expectation that their returns will converge, thereby generating profits. A similar approach can be applied to two stocks with strongly negative correlations.

To quantify the relationship between two stocks, we employ a linear regression model with time-varying intercept and slope: $y_{1, t} = \alpha_t + \beta_t y_{2, t} + \varepsilon_t$. Since the intercept $\alpha_t$ and the slope $\beta_t$ are unknown, we treat them as unobservable state variables and assume a parsimonious dynamic:
\begin{equation}
	\alpha_t = \alpha_{t-1} + \eta_{1, t}, \quad \beta_t = \beta_{t-1} + \eta_{2,t},
\end{equation}
where $\eta_{1, t}$ and $\eta_{2,t}$ are process noises.

Using the language of Kalman filters, the state matrix is $A_t = I_2$ and the observation matrix $C_t = (1, y_{2,t})$. Since the actual noise covariance matrices are unknown, we assume the reference matrices as $Q_t = I_2$ and $R_t = 1$ in the absence of prior information. However, this assumption inevitably results in misspecification.

We analyze the daily closing prices of Google (Alphabet) and Amazon from August 9, 2019, to January 26, 2023, encompassing a total of 873 trading days. Using the first 100 days of data in 2019, we estimate the initial intercept $\alpha_0$ and slope $\beta_0$. The initial wealth is $\$10,000$. We exploit the mean-reverting behavior in the spread $y_{1, t} - \alpha_t - \beta_t y_{2, t}$. When the spread deviates more than 2.0 standard deviations from its mean, we buy or sell $100 \times (1, -\beta_{t})$ shares. The mean and STD are calculated using a rolling window method with a window size of 20. Consequently, trading can commence on the 21st trading day of 2021. We assume a transaction cost of 0.01\% of the value traded.

Table \ref{tab:nonrob_ratio} displays the Sharpe ratio, Sortino ratio, and terminal wealth for individual stocks and the non-robust strategy, which directly employs the Kalman filter. We assume a 2\% annual risk-free interest rate during this period. Table \ref{tab:robust_ratio} presents the results of the three methods. The KL method exhibits marginal improvement over the non-robust counterpart. Conversely, both the OT and BCOT strategies enhance portfolio performance when an appropriate radius is employed. Notably, among the ten radius cases examined, the BCOT method attains the highest Sharpe ratio of $1.19$, Sortino ratio of $3.04$, and terminal wealth of $\$17,533$ when the radius is set at $1.0$. These values exceed the results of both KL and OT across all ten cases.

\begin{table}
	\footnotesize
	\centering
	\begin{tabular}{cccc}
		\hline
		& Non-robust   &  GOOG  &  AMZN \\
		\hline
		Sharpe & 0.9090 &  0.4671 & 0.1813 \\
		\hline
		Sortino & 2.2069 & 0.6545  & 0.2657 \\
		\hline
		Terminal wealth & 15988 & 14504 & 10455 \\
		\hline
	\end{tabular}
	\caption{Terminal wealth and annual Sharpe and Sortino ratios for stocks and the non-robust strategy}\label{tab:nonrob_ratio}
\end{table}

\begin{table*}
	\scriptsize
	\centering
	\begin{tabular}{ccccccccccc}
		\hline
		Radius & 0.1 & 0.2 & 0.3 & 0.4 & 0.5 & 0.6 & 0.7 & 0.8 & 0.9 & 1.0 \\
		\hline 
		BCOT Sharpe & 0.8432 & 0.8953 & 0.9494 & 0.8973 & 1.0968 & 0.8981 & 1.07  & 1.0703 &   1.1322 & 1.1905 \\
		\hline
		OT Sharpe & 0.9062 & 0.9  & 0.8408 & 0.484 & 0.599 & 0.8284 &  0.8372 &  0.9145 & 1.0568 & 0.6908 \\
		\hline
		KL Sharpe & 0.909 & 0.9091 & 0.8972 & 0.8972 & 0.909 & 0.897 &  0.8969 & 0.8968 & 0.8967 & 0.8967 \\
		\hline 
		BCOT Sortino & 2.0257 & 2.0982 & 2.3016 & 2.1601 & 2.7703 & 2.1588 & 2.863 & 2.7822 & 2.8389 & 3.0438 \\
		\hline
		OT Sortino & 2.1341 & 2.1831 &  2.1723 & 1.0761 & 1.2958 & 1.9475 & 1.9682 & 2.2369 & 2.4592 & 1.5942 \\
		\hline
		KL Sortino & 2.2071 & 2.2065 & 2.1591 & 2.1591 & 2.2054 & 2.1581 & 2.1577 & 2.1572 & 2.1567 & 2.1564 \\
		\hline 
		BCOT terminal wealth & 15588 & 15933 & 16248 & 15914 & 17046 & 15900 & 17049 & 16772 & 17278 & 17533 \\
		\hline 
		OT terminal wealth & 16063 & 15921 & 15253 & 13436 & 13962 & 15850 & 15623 & 16001 & 16868 & 15098 \\
		\hline 
		KL terminal wealth & 15988 & 15989 & 15914 & 15914 & 15990 & 15914 & 15913 & 15913 & 15912 & 15913 \\
		\hline
	\end{tabular}
	\caption{Terminal wealth and annual Sharpe and Sortino ratios for the KL, OT, and BCOT methods}\label{tab:robust_ratio}
\end{table*}

%\bibliographystyle{apalike}
%%\bibliographystyle{IEEEtran}
%\bibliography{ref.bib}

\begin{thebibliography}{}
	
	\bibitem[Anderson and Moore, 1979]{anderson1979}
	Anderson, B.~D. and Moore, J.~B. (1979).
	\newblock {\em Optimal Filtering}.
	\newblock Prentice-Hall, Inc., Englewood Cliffs, N.J.
	
	\bibitem[Backhoff-Veraguas et~al., 2017]{backhoff2017causal}
	Backhoff-Veraguas, J., Beiglbock, M., Lin, Y., and Zalashko, A. (2017).
	\newblock Causal transport in discrete time and applications.
	\newblock {\em SIAM Journal on Optimization}, 27(4):2528--2562.
	
	\bibitem[Benson and Vanderbei, 2003]{benson2003solving}
	Benson, H.~Y. and Vanderbei, R.~J. (2003).
	\newblock Solving problems with semidefinite and related constraints using
	interior-point methods for nonlinear programming.
	\newblock {\em Mathematical Programming}, 95(2):279--302.
	
	\bibitem[Bures, 1969]{bures1969extension}
	Bures, D. (1969).
	\newblock An extension of {Kakutani's} theorem on infinite product measures to
	the tensor product of semifinite $w^*$-algebras.
	\newblock {\em Transactions of the American Mathematical Society},
	135:199--212.
	
	\bibitem[Durovic and Kovacevic, 1999]{durovic1999robust}
	Durovic, Z.~M. and Kovacevic, B.~D. (1999).
	\newblock Robust estimation with unknown noise statistics.
	\newblock {\em IEEE Transactions on Automatic Control}, 44(6):1292--1296.
	
	\bibitem[Gandhi and Mili, 2009]{gandhi2009robust}
	Gandhi, M.~A. and Mili, L. (2009).
	\newblock Robust {Kalman} filter based on a generalized maximum-likelihood-type
	estimator.
	\newblock {\em IEEE Transactions on Signal Processing}, 58(5):2509--2520.
	
	\bibitem[Givens and Shortt, 1984]{givens1984class}
	Givens, C.~R. and Shortt, R.~M. (1984).
	\newblock A class of {Wasserstein} metrics for probability distributions.
	\newblock {\em Michigan Mathematical Journal}, 31(2):231--240.
	
	\bibitem[Hassibi et~al., 1999]{hassibi1999indefinite}
	Hassibi, B., Sayed, A.~H., and Kailath, T. (1999).
	\newblock {\em Indefinite-quadratic estimation and control: A unified approach
		to $H^2$ and $H^\infty$ theories}.
	\newblock SIAM.
	
	\bibitem[Huang et~al., 2017]{huang2017novel}
	Huang, Y., Zhang, Y., Wu, Z., Li, N., and Chambers, J. (2017).
	\newblock A novel adaptive {Kalman} filter with inaccurate process and
	measurement noise covariance matrices.
	\newblock {\em IEEE Transactions on Automatic Control}, 63(2):594--601.
	
	\bibitem[Huang et~al., 2020]{huang2020novel}
	Huang, Y., Zhang, Y., Zhao, Y., Shi, P., and Chambers, J.~A. (2020).
	\newblock A novel outlier-robust {Kalman} filtering framework based on
	statistical similarity measure.
	\newblock {\em IEEE Transactions on Automatic Control}, 66(6):2677--2692.
	
	\bibitem[Idan and Speyer, 2010]{idan2010cauchy}
	Idan, M. and Speyer, J.~L. (2010).
	\newblock Cauchy estimation for linear scalar systems.
	\newblock {\em IEEE Transactions on Automatic Control}, 55(6):1329--1342.
	
	\bibitem[Idan and Speyer, 2013]{idan2013multivariate}
	Idan, M. and Speyer, J.~L. (2013).
	\newblock Multivariate {Cauchy} estimator with scalar measurement and process
	noises.
	\newblock In {\em 52nd IEEE Conference on Decision and Control}, pages
	5016--5023. IEEE.
	
	\bibitem[Lassalle, 2013]{lassalle2013causal}
	Lassalle, R. (2013).
	\newblock Causal transference plans and their {Monge-Kantorovich} problems.
	\newblock {\em arXiv preprint arXiv:1303.6925}.
	
	\bibitem[Levy and Nikoukhah, 2012]{levy2012robust}
	Levy, B.~C. and Nikoukhah, R. (2012).
	\newblock Robust state space filtering under incremental model perturbations
	subject to a relative entropy tolerance.
	\newblock {\em IEEE Transactions on Automatic Control}, 58(3):682--695.
	
	\bibitem[Petersen and Savkin, 1999]{petersen1999robust}
	Petersen, I.~R. and Savkin, A.~V. (1999).
	\newblock {\em Robust {Kalman} filtering for signals and systems with large
		uncertainties}.
	\newblock Springer Science \& Business Media.
	
	\bibitem[Petersen and Pedersen, 2012]{petersen2012matrix}
	Petersen, K.~B. and Pedersen, M.~S. (November 15, 2012).
	\newblock {\em The matrix cookbook}.
	
	\bibitem[Shafieezadeh~Abadeh et~al., 2018]{nips2018kalman}
	Shafieezadeh~Abadeh, S., Nguyen, V.~A., Kuhn, D., and Mohajerin~Esfahani, P.~M.
	(2018).
	\newblock Wasserstein distributionally robust {Kalman} filtering.
	\newblock {\em Advances in Neural Information Processing Systems}, 31.
	
	\bibitem[Simon, 2006]{simon2006optimal}
	Simon, D. (2006).
	\newblock {\em Optimal state estimation: {Kalman, H infinity}, and nonlinear
		approaches}.
	\newblock John Wiley \& Sons.
	
	\bibitem[Speyer et~al., 1974]{speyer1974optimization}
	Speyer, J., Deyst, J., and Jacobson, D. (1974).
	\newblock Optimization of stochastic linear systems with additive measurement
	and process noise using exponential performance criteria.
	\newblock {\em IEEE Transactions on Automatic Control}, 19(4):358--366.
	
	\bibitem[Taghvaei and Mehta, 2020]{taghvaei2020optimal}
	Taghvaei, A. and Mehta, P.~G. (2020).
	\newblock An optimal transport formulation of the ensemble {Kalman} filter.
	\newblock {\em IEEE Transactions on Automatic Control}, 66(7):3052--3067.
	
	\bibitem[Villani, 2009]{villani2009optimal}
	Villani, C. (2009).
	\newblock {\em Optimal transport: old and new}, volume 338.
	\newblock Springer.
	
	\bibitem[Yang et~al., 2022]{yang2022decision}
	Yang, J., Zhang, L., Chen, N., Gao, R., and Hu, M. (2022).
	\newblock Decision-making with side information: A causal transport robust
	approach.
	\newblock {\em Optimization Online}.
	
	\bibitem[Zhang et~al., 2020]{zhang2020identification}
	Zhang, L., Sidoti, D., Bienkowski, A., Pattipati, K.~R., Bar-Shalom, Y., and
	Kleinman, D.~L. (2020).
	\newblock On the identification of noise covariances and adaptive kalman
	filtering: A new look at a 50 year-old problem.
	\newblock {\em IEEE Access}, 8:59362--59388.
	
	\bibitem[Zorzi, 2016]{zorzi2016robust}
	Zorzi, M. (2016).
	\newblock Robust {Kalman} filtering under model perturbations.
	\newblock {\em IEEE Transactions on Automatic Control}, 62(6):2902--2907.
	
	\bibitem[Zorzi, 2020]{zorzi2020optimal}
	Zorzi, M. (2020).
	\newblock Optimal transport between {Gaussian} stationary processes.
	\newblock {\em IEEE Transactions on Automatic Control}, 66(10):4939--4944.
	
\end{thebibliography}

\appendix
\section{Proofs of results}\label{sec:append}
\subsection{Proof of Corollary \ref{cor:Gaussian}}
\begin{proof}
	Define the centralized random variables as $\eta_1 := \xi_1 - m_1 \sim \tilde \p_1 = N(0, M_1)$ and $\eta_2 := \xi_2 - m_2 \sim \tilde \p_2 = N(0, M_2)$. Then the integrand reduces to
	\begin{align*}
		& \xi_1^\top P \xi_1 + \xi_2^\top S \xi_2 + \xi_1^\top U \xi_2 \\
		& \quad = (\eta_1 + m_1)^\top P (\eta_1 + m_1) + (\eta_2 + m_2)^\top S (\eta_2 + m_2) \\
		& \qquad + (\eta_1 + m_1)^\top U (\eta_2 + m_2) \\
		&\quad  =  \eta_1^\top P \eta_1 + \eta_1^\top P m_1 + m^\top_1 P \eta_1 + m^\top_1 P m_1 \\ 
		& \qquad + \eta_2^\top S \eta_2  + m_2^\top S \eta_2  + \eta_2^\top S m_2  + m_2^\top S m_2\\
		& \qquad + \eta_1^\top U \eta_2 + m_1^\top U \eta_2 + \eta_1^\top U m_2 + m_1^\top U m_2.
	\end{align*}
	We only need to calculate the OT problem with two zero-mean distributions, that is,
	\begin{align*}
		& \inf_{ \pi \in \Pi(\tilde \p_1, \tilde \p_2)} \int (\eta_1^\top P \eta_1 + \eta_2^\top S \eta_2 + \eta_1^\top U \eta_2) d\pi = \tr[PM_1] + \tr[SM_2] + \inf_{ \pi \in \Pi(\tilde \p_1, \tilde \p_2)} \int \eta_1^\top U \eta_2 d\pi. 
	\end{align*} 
	
	For the last minimization problem, we can regard $U \eta_2$ as a new normal random variable and apply \cite{givens1984class} to $\eta_1$ and  $U \eta_2$ to show that
	\begin{equation*}
		\inf_{ \pi \in \Pi(\tilde \p_1, \tilde \p_2)} \int \eta_1^\top U \eta_2 d\pi = - \tr \left[\sqrt{\sqrt{M_1} U M_2 U \sqrt{M_1}} \right]. 
	\end{equation*}
\end{proof} 

\subsection{Proof of Corollary \ref{cor:Wasser}}
\begin{proof}
	By \citet[Proposition 5.2]{backhoff2017causal}, the bicausal Wasserstein distance \eqref{eq:bicau-dist} can be solved using the dynamic programming principle. The boundary condition implies the value function in Step 2 is
	\begin{equation*}
		v_2 = c = \|x_{t-1} - \bar{x}_{t-1}\|^2 + \| x_t - \bar{x}_t \|^2 + \| y_t - \bar{y}_t \|^2.
	\end{equation*}
	In Step 1, we obtain the value function
	\begin{align*}
		v_1 & = \inf_{ \substack{\pi_t \in \Pi(\p(x_t, y_t|x_{t-1}, y_{1:t-1}), \\ \qquad \barp(\bar{x}_t, \bar{y}_t |\bar{x}_{t-1}, y_{1:t-1} ))}}  \int v_2 d\pi_t \\
		& = \|x_{t-1} - \bar{x}_{t-1}\|^2 + \inf_{ \substack{ \pi_t \in \Pi( \p(x_t, y_t|x_{t-1}), \\ \qquad \barp(\bar{x}_t, \bar{y}_t |\bar{x}_{t-1}) )}}  \int \left( \| x_t - \bar{x}_t \|^2 + \| y_t - \bar{y}_t \|^2 \right) d\pi_t,
	\end{align*}
	by the Markov property of the reference and alternative models, i.e., 
	\begin{equation*}
		\p(x_t, y_t|x_{t-1}, y_{1:t-1}) = \p(x_t, y_t|x_{t-1}) \text{ and } \barp(\bar{x}_t, \bar{y}_t |\bar{x}_{t-1}, y_{1:t-1}) = \barp(\bar{x}_t, \bar{y}_t |\bar{x}_{t-1}).
	\end{equation*}
	By Corollary \ref{cor:Gaussian}, with two Gaussian distributions $\xi_1 = (x^\top_t, y^\top_t)^\top$ and $\xi_2 = (\bar{x}^\top_t, \bar{y}^\top_t)^\top$, and coefficients $P = I_{n+m}$, $S = I_{n+m}$, $U = - 2I_{n+m}$, we obtain 
	\begin{align*}
		v_1 =&  \|x_{t-1} - \bar{x}_{t-1}\|^2 + \begin{pmatrix}
			x^\top_{t-1} A^\top_t,  x^\top_{t-1} A^\top_t C^\top_t 
		\end{pmatrix}  \begin{pmatrix}
			A_t x_{t-1} \\
			C_t A_t x_{t-1}
		\end{pmatrix} \\
		&+ \begin{pmatrix}
			\bar{x}^\top_{t-1} A^\top_t,  \bar{x}^\top_{t-1} A^\top_t C^\top_t 
		\end{pmatrix}  \begin{pmatrix}
			A_t \bar{x}_{t-1} \\
			C_t A_t \bar{x}_{t-1}
		\end{pmatrix} \\
		& - 2 \begin{pmatrix}
			x^\top_{t-1} A^\top_t,  x^\top_{t-1} A^\top_t C^\top_t 
		\end{pmatrix}  \begin{pmatrix}
			A_t \bar{x}_{t-1} \\
			C_t A_t \bar{x}_{t-1} 
		\end{pmatrix} \\
		& + \tr \Big[\sigma^2(Q_t, R_t) + \sigma^2(\bar{Q}_t, \bar{R}_t) - 2 \sqrt{\sqrt{\sigma^2(Q_t, R_t)} \sigma^2(\bar{Q}_t, \bar{R}_t) \sqrt{\sigma^2(Q_t, R_t)} } \Big].
	\end{align*}
	
	With the definition of $H$, applying Corollary \ref{cor:Gaussian} again yields
	\begin{align*}
		& \cW_{bc} (\p_{t-1}, \barp_{t-1}) = v_0 \\
		& \quad = \inf_{ \pi_{t-1} \in \Pi(\p(x_{t-1}|y_{1:t-1}), \barp(\bar{x}_{t-1} |y_{1:t-1}))}  \int  v_1 d\pi_{t-1} \\
		& \quad  = \hat{x}^\top_{t-1} H \hat{x}_{t-1} + \hat{x}^\top_{t-1} H \hat{x}_{t-1} - 2 \hat{x}^\top_{t-1} H \hat{x}_{t-1} \\
		& \qquad + \tr \Big[\sigma^2(Q_t, R_t) + \sigma^2(\bar{Q}_t, \bar{R}_t) - 2 \sqrt{\sqrt{\sigma^2(Q_t, R_t)} \sigma^2(\bar{Q}_t, \bar{R}_t) \sqrt{\sigma^2(Q_t, R_t)} } \Big] \\
		& \qquad  + \tr \Big[ H \Sigma_{t-1} + H \bar{\Sigma}_{t-1} - 2 \sqrt{\sqrt{\Sigma_{t-1}} H \bar{\Sigma}_{t-1} H \sqrt{\Sigma_{t-1}} } \Big] \\
		& \quad =: w(\bar{Q}_t, \bar{R}_t, \bar{\Sigma}_{t-1}).
	\end{align*}
\end{proof}

\subsection{Proof of Theorem \ref{thm:covex}}
\begin{proof}
	First, we prove the infimum and supremum in the problem \eqref{eq:minimax} are interchangeable. Moreover, the infimum over $f \in \cL$ can be restricted to linear functions.
	
	By the weak duality and the fact that the optimizer is linear when $\barp$ is fixed, we obtain
	\begin{align*}
		& \inf_{f \in \cL} \sup_{\barp \in \cB_{\varepsilon, t - 1}} \E_{\barp} \left[ \| \bar{x}_t - f(\bar{y}_t) \|^2 \Big| y_{1:t-1} \right] \\
		& \quad \geq \sup_{\barp \in \cB_{\varepsilon, t - 1}} \inf_{f \in \cL} \E_{\barp} \left[ \| \bar{x}_t - f(\bar{y}_t) \|^2 \Big| y_{1:t-1} \right] \\
		& \quad =  \sup_{\barp \in \cB_{\varepsilon, t - 1}} \inf_{G, g} \E_{\barp} \left[ \| \bar{x}_t - G \bar{y}_t - g \|^2 \Big| y_{1:t-1} \right].
	\end{align*}
	Moreover, since a linear estimator is in $\cL$, then
	\begin{align*}
		& \inf_{f \in \cL} \sup_{\barp \in \cB_{\varepsilon, t - 1}} \E_{\barp} \left[ \| \bar{x}_t - f(\bar{y}_t) \|^2 \Big| y_{1:t-1} \right] \leq \inf_{G, g} \sup_{\barp \in \cB_{\varepsilon, t - 1}}  \E_{\barp} \left[ \| \bar{x}_t - G \bar{y}_t - g \|^2 \Big| y_{1:t-1} \right].
	\end{align*}
	If we can prove the strong duality for linear estimators, that is,
	\begin{align}
		& \sup_{\barp \in \cB_{\varepsilon, t - 1}}  \inf_{G, g} \E_{\barp} \left[ \| \bar{x}_t - G \bar{y}_t - g \|^2 \Big| y_{1:t-1} \right]  = \inf_{G, g} \sup_{\barp \in \cB_{\varepsilon, t - 1}}  \E_{\barp} \left[ \| \bar{x}_t - G \bar{y}_t - g \|^2 \Big| y_{1:t-1} \right], \label{eq:strong-dual-linear}
	\end{align}
	then all inequalities above are equalities and in particular,
	\begin{align}
		& \inf_{f \in \cL} \sup_{\barp \in \cB_{\varepsilon, t - 1}} \E_{\barp} \left[ \| \bar{x}_t - f(\bar{y}_t) \|^2 \Big| y_{1:t-1} \right] =	\sup_{\barp \in \cB_{\varepsilon, t - 1}}  \inf_{G, g} \E_{\barp} \left[ \| \bar{x}_t - G \bar{y}_t - g \|^2 \Big| y_{1:t-1} \right].  \label{eq:exchanged}  
	\end{align}
	To prove the strong duality \eqref{eq:strong-dual-linear}, we invoke Sion's minimax theorem. Note that $\barp$ is fully characterized by the choice of $\bar{Q}_t, \bar{R}_t, \bar{\Sigma}_{t-1}$. With a given $\barp$, we have
	\begin{align*}
		\E_{\barp} [ \bar{x}^\top_t \bar{x}_t | y_{1:t-1} ] = & \hat{x}^\top_{t-1} A^\top_t A_t \hat{x}_{t-1} + \tr[A_t \bar{\Sigma}_{t-1} A^\top_t + \bar{Q}_t], \\
		\E_{\barp} [ \bar{y}^\top_t G^\top G \bar{y}_t | y_{1:t-1} ] = & \hat{x}^\top_{t-1} A^\top_t C^\top_t G^\top G C_t A_t \hat{x}_{t-1} + \tr[ C_t A_t \bar{\Sigma}_{t-1} A^\top_t C^\top_t G^\top G] \\
		& + \tr[(C_t \bar{Q}_t C^\top_{t} + \bar{R}_t) G^\top G], \\
		\E_{\barp} [ \bar{y}^\top_t G^\top \bar{x}_t | y_{1:t-1} ] = & \E_{\barp} [ \bar{x}^\top_t G \bar{y}_t | y_{1:t-1} ] = \hat{x}^\top_{t-1} A^\top_t G C_t A_t \hat{x}_{t-1} + \tr[ (C_t A_t \bar{\Sigma}_{t-1} A^\top_t + C_t \bar{Q}_t) G], \\
		\E_{\barp} [ \bar{y}^\top_t G^\top g | y_{1:t-1} ] = & \E_{\barp} [ g^\top G \bar{y}_t | y_{1:t-1} ]  =  \hat{x}^\top_{t-1} A^\top_t C^\top_t G^\top g, \\
		\E_{\barp} [ g^\top \bar{x}_t | y_{1:t-1} ] = & \E_{\barp} [ \bar{x}^\top_t g| y_{1:t-1} ] = \hat{x}^\top_{t-1} A^\top_t g.
	\end{align*} 
	Then we verify the conditions for Sion's minimax theorem.
	\begin{enumerate}
		\item When written in a matrix form, the objective is linear and thus quasi-concave in $\bar{Q}_t, \bar{R}_t, \bar{\Sigma}_{t-1}$. Besides, it is convex and thus quasi-convex in $G, g$. 
		\item The domain for $G, g$ is convex.
		\item Since the set $\{\bar{Q}_t \succeq 0, \bar{R}_t \succeq \delta I_m, \bar{\Sigma}_{t-1} \succeq 0 \}$ is closed and convex, we only need to prove the set $\{ w(\bar{Q}_t, \bar{R}_t, \bar{\Sigma}_{t-1}) \leq \varepsilon \}$ is compact and convex. Noting that the trace operator is linear and $\sigma^2(\bar{Q}_t, \bar{R}_t) $ is linear in $(\bar{Q}_t, \bar{R}_t)$, it suffices to show that the function $\tr[\sqrt{(\cdot)}]$ is concave. Consider two positive semi-definite matrices $X$ and $Y$. We want to show 
		\begin{equation*}
			\tr[\sqrt{\lambda X + (1 - \lambda) Y}] \geq \lambda \tr[\sqrt{X}] +  (1 - \lambda) \tr[\sqrt{Y}].
		\end{equation*}
		Equivalently, we need to prove
		\begin{equation*}
			\sqrt{\lambda X + (1 - \lambda) Y} \succeq \lambda \sqrt{X} +  (1 - \lambda) \sqrt{Y}.
		\end{equation*}
		Indeed, the square root is monotone for positive semi-definite matrices. Moreover,
		\begin{align*}
			& \lambda X + (1 - \lambda) Y - (\lambda \sqrt{X} +  (1 - \lambda) \sqrt{Y}) (\lambda \sqrt{X} +  (1 - \lambda) \sqrt{Y}) \\
			& = \lambda (1 - \lambda) (X + Y - \sqrt{X} \sqrt{Y} - \sqrt{Y} \sqrt{X}) \succeq 0.
		\end{align*}
		Then $w(\bar{Q}_t, \bar{R}_t, \bar{\Sigma}_{t-1})$ is convex and $w(\bar{Q}_t, \bar{R}_t, \bar{\Sigma}_{t-1}) \leq \varepsilon$ is a convex set. It is straightforward to see the set is compact. 
	\end{enumerate}
	Conditions in Sion's minimax theorem are verified. Hence, \eqref{eq:strong-dual-linear} and \eqref{eq:exchanged} hold.
	
	Next, we further simplify the right-hand side of \eqref{eq:exchanged}. Minimizing over $g$ and setting the first derivative to be zero, we obtain
	\begin{equation*}
		g_* = A_t \hat{x}_{t-1} - G C_t A_t \hat{x}_{t-1}.
	\end{equation*}
	On the right-hand side of \eqref{eq:exchanged}, terms depending on $g_*$ become
	\begin{align*}
		& (g_*)^\top g_* + 2 \hat{x}^\top_{t-1} A^\top_t C^\top_t G^\top g_* - 2 \hat{x}^\top_{t-1} A^\top_t g_* \\
		& \quad = -  \hat{x}^\top_{t-1} A^\top_t A_t \hat{x}_{t-1} + \hat{x}^\top_{t-1} A^\top_t C^\top_t G^\top A_t \hat{x}_{t-1} \\
		& \qquad + \hat{x}^\top_{t-1} A^\top_t G C_t A_t \hat{x}_{t-1}  -  \hat{x}^\top_{t-1} A^\top_t C^\top_t G^\top G C_t A_t \hat{x}_{t-1}. 
	\end{align*}
	The objective $\E_{\barp} \left[ \| \bar{x}_t - G \bar{y}_t - g \|^2 \Big| y_{1:t-1} \right]$ with $g_*$ is given by
	\begin{align*}
		& \tr[A_t \bar{\Sigma}_{t-1} A^\top_t + \bar{Q}_t] + \tr[ C_t A_t \bar{\Sigma}_{t-1} A^\top_t C^\top_t G^\top G] \\
		& + \tr[(C_t \bar{Q}_t C^\top_{t} + \bar{R}_t) G^\top G] - \tr[ (C_t A_t \bar{\Sigma}_{t-1} A^\top_t + C_t \bar{Q}_t) G] \\ 
		&- \tr[ G^\top (A_t\bar{\Sigma}_{t-1} A^\top_t C^\top_t + \bar{Q}_t C^\top_t)].
	\end{align*}
	Or in a matrix form, it becomes
	\begin{equation}\label{eq:G-obj}
		\tr \left[
		\begin{pmatrix}
			E_{11} & E_{12} \\
			E^\top_{12} & E_{22} 
		\end{pmatrix} \begin{pmatrix}
			I_n, & -G \\
			-G^\top, & G^\top G
		\end{pmatrix} \right],
	\end{equation}
	with
	\begin{align*}
		E_{11}  = & A_t \bar{\Sigma}_{t-1} A^\top_t + \bar{Q}_t, \\
		E_{12} = & A_t\bar{\Sigma}_{t-1} A^\top_t C^\top_t + \bar{Q}_t C^\top_t, \\
		E_{22} = & C_t A_t \bar{\Sigma}_{t-1} A^\top_t C^\top_t + C_t \bar{Q}_t C^\top_{t} + \bar{R}_t.
	\end{align*}
	
	To minimize over $G$, note that $\bar{R}_{t} \succeq \delta I_m$ with a given sufficiently small $\delta > 0$. Thus, the objective \eqref{eq:G-obj} is strictly convex in $G$. The minimization problem over $G$ has a unique solution $G_*$, which can be solved from the first-order optimality condition:
	\begin{equation*}
		2 G_* K - 2 M^\top = 0. 
	\end{equation*} 
	Therefore, $G_* = M^\top K^{-1}$. Substituting $G_*$ into the objective \eqref{eq:G-obj}, we obtain \eqref{eq:main-pro}. The robust optimizer \eqref{rob_f} is proved with the expressions of $G_*$ and $g_*$. 
\end{proof} 

\subsection{Proof of Lemma \ref{lem:convex}}
\begin{proof}
	The idea is to investigate the Hessian matrix of the objective function $F$ in \eqref{eq:main-pro}. Since $\tr[ A_t \bar{\Sigma}_{t-1} A^\top_t + \bar{Q}_t] $ is linear, we only need to show the Hessian of $\tr[ M M^\top K^{-1}]$ is positive semi-definite. We omit the time subscript for simplicity.

	{\noindent \bf Derivatives on $\bar{Q}$} 
	
	As $\tr[M^\top K^{-1} M] = \sum_r (M^\top K^{-1} M)_{rr}$, we have
	\begin{align*}
		& \frac{\partial}{\partial \bar{Q}_{ij}} (M^\top K^{-1} M)_{rr} = \frac{\partial}{\partial \bar{Q}_{ij}} \sum_{w, x} (M^\top)_{rw} (K^{-1})_{wx} M_{xr} \\
		& =  \sum_{w, x} \Big\{ \frac{\partial (M^\top)_{rw}}{\partial \bar{Q}_{ij}}   (K^{-1})_{wx} M_{xr} + (M^\top)_{rw} \frac{\partial (K^{-1})_{wx}}{\partial \bar{Q}_{ij}} M_{xr} + (M^\top)_{rw} (K^{-1})_{wx} \frac{ \partial M_{xr}}{\partial \bar{Q}_{ij}} \Big\}.
	\end{align*}
	Since $M_{xr} = \text{constant} + \sum_z C_{xz} \bar{Q}_{zr}$, then
	\begin{align*}
		\frac{ \partial M_{xr}}{\partial \bar{Q}_{ij}} = C_{xi} \delta_{rj}, \quad \frac{ \partial (M^\top)_{rw}}{\partial \bar{Q}_{ij}} = C_{wi} \delta_{rj},
	\end{align*}
	where $\delta_{rj}$ is the Kronecker delta.
	
	By the formula for the derivative of matrix inverse, see \citet[Equation 60]{petersen2012matrix}, we have
	\begin{align*}
		\frac{\partial (K^{-1})_{wx}}{\partial \bar{Q}_{ij}} & = \sum_{v, z} \frac{\partial (K^{-1})_{wx}}{\partial K_{vz}} \frac{\partial K_{vz}}{\partial \bar{Q}_{ij}}  = - \sum_{v, z} (K^{-1})_{wv} (K^{-1})_{zx} C_{vi} (C^\top)_{jz} = - (K^{-1} C)_{wi} (C^\top K^{-1})_{jx}.
	\end{align*}
	Then
	\begin{align*}
		& \frac{\partial}{\partial \bar{Q}_{ij}} (M^\top K^{-1} M)_{rr} \\
		& = \sum_{w, x} C_{wi} \delta_{rj}  (K^{-1})_{wx} M_{xr} \\
		& \quad - \sum_{v, z, w, x} (M^\top)_{rw}  (K^{-1})_{wv} (K^{-1})_{zx} C_{vi} (C^\top)_{jz} M_{xr} +  \sum_{w, x}  (M^\top)_{rw} (K^{-1})_{wx} C_{xi} \delta_{rj}.
	\end{align*}
	Next, we sum over $r$ and rearrange terms in an order of matrix multiplication. Since $K^{-1}$ is symmetric, we obtain
	\begin{align*}
		& \frac{\partial}{\partial \bar{Q}_{ij}} \tr[M^\top K^{-1} M] \\
		& = 2 \sum_{w, x} (C^\top)_{iw} (K^{-1})_{wx} M_{xj}  - \sum_{r, v, z, w, x} (C^\top)_{iv} (K^{-1})_{vw} M_{wr}  (M^\top)_{rx} (K^{-1})_{xz} C_{zj}.
	\end{align*}
	For the second derivatives with respect to $\bar{Q}$, we have
	\begin{align*}
		& \frac{\partial^2}{\partial \bar{Q}_{ij} \partial \bar{Q}_{kl}} \tr[M^\top K^{-1} M] \\
		& = - 2 \sum_{w, x} (C^\top)_{iw} (K^{-1} C)_{wk} (C^\top K^{-1})_{lx} M_{xj} \\
		& \quad + 2 \sum_{w, x} (C^\top)_{iw} (K^{-1})_{wx} C_{xk} \delta_{lj} \\
		& \quad - \sum_{r, v, z, w, x} (C^\top)_{iv} (-1) (K^{-1} C)_{vk} (C^\top K^{-1})_{lw} M_{wr}  \\
		& \qquad \qquad \qquad  \times (M^\top)_{rx} (K^{-1})_{xz} C_{zj} \\
		&\quad  - \sum_{r, v, z, w, x} (C^\top)_{iv} (K^{-1})_{vw} C_{wk} \delta_{lr}  (M^\top)_{rx} (K^{-1})_{xz} C_{zj} \\
		& \quad - \sum_{r, v, z, w, x} (C^\top)_{iv} (K^{-1})_{vw} M_{wr}  C_{xk} \delta_{lr} (K^{-1})_{xz} C_{zj} \\
		& \quad  - \sum_{r, v, z, w, x} (C^\top)_{iv} (K^{-1})_{vw} M_{wr}  (M^\top)_{rx} \\
		& \qquad \qquad \qquad \times (-1) (K^{-1} C)_{xk} (C^\top K^{-1})_{lz} C_{zj} \\
		& = - 2 (C^\top K^{-1} C)_{ik} (C^\top K^{-1} M)_{lj} + 2 (C^\top K^{-1} C)_{ik} \delta_{lj} \\
		& \quad + (C^\top K^{-1} C)_{ik} (C^\top K^{-1} M M^\top K^{-1} C)_{lj} \\
		& \quad - (C^\top K^{-1} C)_{ik} (M^\top K^{-1} C)_{lj} \\
		& \quad - (C^\top K^{-1} M)_{il} (C^\top K^{-1} C)_{kj} \\
		& \quad + (C^\top K^{-1} M M^\top K^{-1} C)_{ik} (C^\top K^{-1} C)_{lj}.
	\end{align*}
	Besides,
	\begin{align*}
		& \frac{\partial^2}{\partial \bar{Q}_{ij} \partial \bar{\Sigma}_{kl}} \tr[M^\top K^{-1} M] \\
		& = - 2 \sum_{w, x} (C^\top)_{iw} (K^{-1} C A)_{wk} (A^\top C^\top K^{-1})_{lx} M_{xj} \\
		& \quad + 2 \sum_{w, x} (C^\top)_{iw} (K^{-1})_{wx} (CA)_{xk} (A^\top)_{lj} \\
		& \quad - \sum_{r, v, z, w, x} (C^\top)_{iv} (-1) (K^{-1} C A)_{vk} ( A^\top C^\top K^{-1})_{lw} \\
		& \qquad \qquad \qquad \times M_{wr}  (M^\top)_{rx} (K^{-1})_{xz} C_{zj} \\
		& \quad - \sum_{r, v, z, w, x} (C^\top)_{iv} (K^{-1})_{vw} (CA)_{wk} \\
		& \qquad \qquad \qquad \times (A^\top)_{lr}  (M^\top)_{rx} (K^{-1})_{xz} C_{zj} \\
		& \quad - \sum_{r, v, z, w, x} (C^\top)_{iv} (K^{-1})_{vw} M_{wr} (CA)_{xk} \\
		& \qquad \qquad \qquad \times (A^\top)_{lr} (K^{-1})_{xz} C_{zj} \\
		& \quad - \sum_{r, v, z, w, x} (C^\top)_{iv} (K^{-1})_{vw} M_{wr}  (M^\top)_{rx} \\
		& \qquad \qquad \qquad \times (-1) (K^{-1} C A)_{xk} (A^\top C^\top K^{-1})_{lz} C_{zj} \\
		& =  - 2 (C^\top K^{-1} C A)_{ik} (A^\top C^\top K^{-1} M)_{lj} \\
		& \quad + 2 (C^\top K^{-1} C A)_{ik} (A^\top)_{lj} \\
		& \quad + (C^\top K^{-1} C A)_{ik} (A^\top C^\top K^{-1} M M^\top K^{-1} C)_{lj} \\
		& \quad - (C^\top K^{-1} C A)_{ik} (A^\top M^\top K^{-1} C)_{lj} \\
		& \quad - (C^\top K^{-1} M A)_{il} (A^\top C^\top K^{-1} C)_{kj} \\
		& \quad + (C^\top K^{-1} M M^\top K^{-1} C A)_{ik} (A^\top C^\top K^{-1} C)_{lj}.
	\end{align*}
	
	{\noindent \bf Derivatives on $\bar{R}$} 
	
	Note that $\bar{R}$ appears only in $K$. The first derivative is
	\begin{align*}
		& \frac{\partial}{\partial \bar{R}_{ij}} \tr[M^\top K^{-1} M] = - \sum_{w, x, y} (K^{-1})_{iw} M_{wx} (M^\top)_{xy} (K^{-1})_{yj}.
	\end{align*}
	The second derivatives are
	\begin{align*}
		& \frac{\partial^2}{\partial \bar{R}_{ij} \partial \bar{R}_{kl}} \tr[M^\top K^{-1} M] \\
		& = -  \sum_{w, x, y} (-1) (K^{-1})_{ik} (K^{-1})_{wl} M_{wx} (M^\top)_{xy} (K^{-1})_{yj} \\
		&\quad - \sum_{w, x, y} (-1) (K^{-1})_{iw} M_{wx} (M^\top)_{xy} (K^{-1})_{yk} (K^{-1})_{lj} \\
		& = (K^{-1})_{ik} (K^{-1} M M^\top K^{-1})_{lj} + (K^{-1} M M^\top K^{-1})_{ik} (K^{-1})_{lj}
	\end{align*}
	and
	\begin{align*}
		&\frac{\partial^2}{\partial \bar{R}_{ij} \partial \bar{Q}_{kl}} \tr[M^\top K^{-1} M] \\
		& = \sum_{w, x, y} (K^{-1} C)_{ik} (C^\top K^{-1})_{lw} M_{wx} (M^\top)_{xy} (K^{-1})_{yj} \\
		& \quad - \sum_{w, x, y} (K^{-1})_{iw} C_{wk} \delta_{lx} (M^\top)_{xy} (K^{-1})_{yj} \\
		&  \quad - \sum_{w, x, y} (K^{-1})_{iw} M_{wx} C_{yk} \delta_{lx} (K^{-1})_{yj} \\
		& \quad + \sum_{w, x, y} (K^{-1})_{iw} M_{wx} (M^\top)_{xy} (K^{-1} C)_{yk} (C^\top K^{-1})_{lj} \\
		& = (K^{-1} C)_{ik} (C^\top K^{-1} M M^\top K^{-1})_{lj} - (K^{-1} C)_{ik} (M^\top K^{-1})_{lj} \\
		& \quad - (C^\top K^{-1})_{kj} (K^{-1} M)_{il} + (K^{-1} MM^\top K^{-1} C)_{ik} (C^\top K^{-1})_{lj}.
	\end{align*}
	With 
	\begin{align*}
		\frac{\partial (K^{-1})_{iw}}{\partial \bar{\Sigma}_{kl}} & = - (K^{-1} C A)_{ik} (A^\top C^\top K^{-1})_{lw}, \\
		\frac{\partial M_{wx}}{\partial \bar{\Sigma}_{kl}} & = (C A)_{wk} (A^\top)_{lx},
	\end{align*}
	we obtain
	\begin{align*}
		& \frac{\partial^2}{\partial \bar{R}_{ij} \partial \bar{\Sigma}_{kl}} \tr[M^\top K^{-1} M] \\
		& = \sum_{w, x, y} (K^{-1} C A)_{ik} (A^\top C^\top K^{-1})_{lw} M_{wx} (M^\top)_{xy} (K^{-1})_{yj} \\
		& \quad - \sum_{w, x, y} (K^{-1})_{iw} (CA)_{wk} (A^\top)_{lx} (M^\top)_{xy} (K^{-1})_{yj} \\
		& \quad - \sum_{w, x, y} (K^{-1})_{iw} M_{wx} (CA)_{yk} (A^\top)_{lx} (K^{-1})_{yj} \\
		&\quad + \sum_{w, x, y} (K^{-1})_{iw} M_{wx} (M^\top)_{xy} (K^{-1} C A)_{yk} (A^\top C^\top K^{-1})_{lj} \\
		& = (K^{-1} C A)_{ik} (A^\top C^\top K^{-1} M M^\top K^{-1})_{lj} \\
		& \quad - (K^{-1} C A)_{ik} (A^\top M^\top K^{-1})_{lj} \\
		& \quad - (A^\top C^\top K^{-1})_{kj} (K^{-1} M A)_{il} \\
		& \quad + (K^{-1} MM^\top K^{-1} C A)_{ik} (A^\top C^\top K^{-1})_{lj}.
	\end{align*}
	
	{\noindent \bf Derivatives on $\bar{\Sigma}$} 
	
	$\bar{\Sigma}$ appears in a similar position of $\bar{Q}$. Therefore,
	\begin{align*}
		\frac{\partial}{\partial \bar{\Sigma}_{ij}} \tr[M^\top K^{-1} M]  = & 2 \sum_{w, x, y} (A^\top C^\top)_{iw} (K^{-1})_{wx} M_{xy} A_{yj}  \\
		&- \sum_{r, v, z, w, x} (A^\top C^\top)_{iv} (K^{-1})_{vw} M_{wr}  (M^\top)_{rx} (K^{-1})_{xz} (CA)_{zj}.
	\end{align*}
	For the second derivatives, we have 
	\begin{align*}
		& \frac{\partial^2}{\partial \bar{\Sigma}_{ij} \partial \bar{\Sigma}_{kl}} \tr[M^\top K^{-1} M] \\
		& = - 2 (A^\top C^\top K^{-1} C A)_{ik} (A^\top C^\top K^{-1} M A)_{lj} \\
		& \quad + 2 (A^\top C^\top K^{-1} C A)_{ik} (A^\top A)_{lj} \\
		& \quad + (A^\top C^\top K^{-1} C A)_{ik} (A^\top C^\top K^{-1} M M^\top K^{-1} C A)_{lj} \\
		& \quad - (A^\top C^\top K^{-1} C A)_{ik} (A^\top M^\top K^{-1} C A)_{lj} \\
		& \quad - (A^\top C^\top K^{-1} M A)_{il} (A^\top C^\top K^{-1} C A)_{kj} \\
		& \quad + (A^\top C^\top K^{-1} M M^\top K^{-1} C A)_{ik} (A^\top C^\top K^{-1} C A)_{lj}.
	\end{align*}
	
	To validate that the Hessian matrix is positive semi-definite, we take arbitrary positive semi-definite matrices $Q, R, \Sigma$ and view them as double-indexed vectors. Consider the quadratic form
	\begin{align}
		\sum_{i, j, k, l} \Big\{ & Q_{ij}  \frac{\partial^2 \tr[M^\top K^{-1} M]}{\partial \bar{Q}_{ij} \partial \bar{Q}_{kl}} Q_{kl} + Q_{ij} \frac{\partial^2 \tr[M^\top K^{-1} M]}{\partial \bar{Q}_{ij} \partial \bar{R}_{kl}} R_{kl} \nonumber \\
		& + Q_{ij} \frac{\partial^2 \tr[M^\top K^{-1} M]}{\partial \bar{Q}_{ij} \partial \bar{\Sigma}_{kl}} \Sigma_{kl} + R_{ij} \frac{\partial^2 \tr[M^\top K^{-1} M]}{\partial \bar{R}_{ij} \partial \bar{R}_{kl}} R_{kl} \nonumber \\
		& + R_{ij} \frac{\partial^2 \tr[M^\top K^{-1} M]}{\partial \bar{R}_{ij} \partial \bar{Q}_{kl}} Q_{kl} + R_{ij} \frac{\partial^2 \tr[M^\top K^{-1} M]}{\partial \bar{R}_{ij} \partial \bar{\Sigma}_{kl}} \Sigma_{kl} \label{eq:Hess} \\
		& + \Sigma_{ij} \frac{\partial^2 \tr[M^\top K^{-1} M]}{\partial \bar{\Sigma}_{ij} \partial \bar{\Sigma}_{kl}} \Sigma_{kl} + \Sigma_{ij} \frac{\partial^2 \tr[M^\top K^{-1} M]}{\partial \bar{\Sigma}_{ij} \partial \bar{Q}_{kl}} Q_{kl} \nonumber \\
		& + \Sigma_{ij} \frac{\partial^2 \tr[M^\top K^{-1} M]}{\partial \bar{\Sigma}_{ij} \partial \bar{R}_{kl}} R_{kl} \Big\}. \nonumber
	\end{align}
	
	We take the first term $ Q_{ij}  \frac{\partial^2 \tr[M^\top K^{-1} M]}{\partial \bar{Q}_{ij} \partial \bar{Q}_{kl}} Q_{kl}$ as an example. We have
	\begin{align*}
		& \sum_{i, j, k, l}  -2 (C^\top K^{-1} C)_{ik} Q_{kl} (C^\top K^{-1} M)_{lj} Q_{ij}  \\
		& = - 2 \tr \left[ C^\top K^{-1} C Q C^\top K^{-1} M Q^\top \right] \\
		& = - 2 \tr \left[ C^\top K^{-1} C Q C^\top K^{-1} M Q\right],
	\end{align*}
	where the second equality follows from that $Q$ is symmetric. Other terms can be calculated in the same way. Then
	\begin{align*}
		& \sum_{i, j, k, l} Q_{ij}  \frac{\partial^2 \tr[M^\top K^{-1} M]}{\partial \bar{Q}_{ij} \partial \bar{Q}_{kl}} Q_{kl} \\
		& = - 2 \tr \left[ C^\top K^{-1} C Q C^\top K^{-1} M Q \right] + 2 \tr \left[ C^\top K^{-1} C Q Q \right] + \tr \left[ C^\top K^{-1} C Q C^\top K^{-1} M M^\top K^{-1} C Q \right] \\
		& \quad - \tr\left[ C^\top K^{-1} C Q M^\top K^{-1} C Q \right] - \tr \left[ C^\top K^{-1} M Q C^\top K^{-1} C Q \right] \\
		&\quad + \tr \left[ C^\top K^{-1} M M^\top K^{-1} C Q C^\top K^{-1} C Q \right].
	\end{align*}
	Similarly, we have the following results for other terms:
	\begin{align*}
		& \sum_{i, j, k, l} R_{ij} \frac{\partial^2 \tr[M^\top K^{-1} M]}{\partial \bar{R}_{ij} \partial \bar{R}_{kl}} R_{kl} \\
		& = \tr[ K^{-1} R K^{-1} M M^\top K^{-1} R + K^{-1} M M^\top K^{-1} R K^{-1} R].
	\end{align*}
	
	\begin{align*}
		& \sum_{i, j, k, l} \Sigma_{ij} \frac{\partial^2 \tr[M^\top K^{-1} M] }{\partial \bar{\Sigma}_{ij} \partial \bar{\Sigma}_{kl}} \Sigma_{kl} \\
		& = - 2 \tr[ A^\top C^\top K^{-1} C A \Sigma A^\top C^\top K^{-1} M A \Sigma] \\
		& \quad + 2 \tr[A^\top C^\top K^{-1} C A \Sigma A^\top A \Sigma] \\
		& \quad + \tr[ A^\top C^\top K^{-1} C A \Sigma A^\top C^\top K^{-1} M M^\top K^{-1} C A \Sigma] \\
		& \quad - \tr[ A^\top C^\top K^{-1} C A \Sigma A^\top M^\top K^{-1} C A \Sigma] \\
		& \quad - \tr [A^\top C^\top K^{-1} M A \Sigma A^\top C^\top K^{-1} C A \Sigma] \\
		& \quad + \tr[ A^\top C^\top K^{-1} M M^\top K^{-1} C A \Sigma A^\top C^\top K^{-1} C A \Sigma].
	\end{align*}
	
	\begin{align*}
		& \sum_{i, j, k, l} \Big\{ R_{ij} \frac{\partial^2 \tr[M^\top K^{-1} M]}{\partial \bar{R}_{ij} \partial \bar{Q}_{kl}} Q_{kl} + Q_{ij} \frac{\partial^2 \tr[M^\top K^{-1} M]}{\partial \bar{Q}_{ij} \partial \bar{R}_{kl}} R_{kl} \Big\}\\
		& = 2 \tr[ K^{-1} C Q C^\top K^{-1} M M^\top K^{-1} R] - 2 \tr[ K^{-1} C Q M^\top K^{-1} R] \\
		& \quad - 2 \tr[ C^\top K^{-1} R K^{-1} M Q] + 2 \tr[ K^{-1} MM^\top K^{-1} C Q C^\top K^{-1} R].
	\end{align*}
	
	\begin{align*}
		& \sum_{i, j, k, l} \Big\{ Q_{ij} \frac{\partial^2 \tr[M^\top K^{-1} M]}{\partial \bar{Q}_{ij} \partial \bar{\Sigma}_{kl}} \Sigma_{kl} + \Sigma_{ij} \frac{\partial^2 \tr[M^\top K^{-1} M]}{\partial \bar{\Sigma}_{ij} \partial \bar{Q}_{kl}} Q_{kl} \Big\} \\
		& = - 4 \tr[C^\top K^{-1} C A \Sigma A^\top C^\top K^{-1} M Q] + 4 \tr[ C^\top K^{-1} C A \Sigma A^\top Q] \\
		& \quad + 2 \tr [ C^\top K^{-1} C A \Sigma A^\top C^\top K^{-1} M M^\top K^{-1} C Q] \\
		& \quad - 2 \tr[ C^\top K^{-1} C A \Sigma A^\top M^\top K^{-1} C Q] \\
		& \quad - 2 \tr[ C^\top K^{-1} M A \Sigma A^\top C^\top K^{-1} C Q] \\
		& \quad + 2 \tr[ C^\top K^{-1} M M^\top K^{-1} C A \Sigma A^\top C^\top K^{-1} C Q].
	\end{align*}
	
	\begin{align*}
		&\sum_{i, j, k, l} \Big\{ R_{ij} \frac{\partial^2 \tr[M^\top K^{-1} M]}{\partial \bar{R}_{ij} \partial \bar{\Sigma}_{kl}} \Sigma_{kl} + \Sigma_{ij} \frac{\partial^2 \tr[M^\top K^{-1} M]}{\partial \bar{\Sigma}_{ij} \partial \bar{R}_{kl}} R_{kl} \Big\} \\
		& = 2 \tr[ K^{-1} C A \Sigma A^\top C^\top K^{-1} M M^\top K^{-1} R] \\
		& \quad - 2 \tr[K^{-1} C A \Sigma A^\top M^\top K^{-1} R] \\
		& \quad - 2 \tr[ K^{-1} M A \Sigma A^\top C^\top K^{-1} R] \\
		& \quad + 2 \tr[ K^{-1} MM^\top K^{-1} C A \Sigma A^\top C^\top K^{-1} R].
	\end{align*}
	Summing all terms up, we find that the quadratic form \eqref{eq:Hess} equals to $2 \tr[K^{-1} L L^\top]$ with
	\begin{align*}	
		L := & R K^{-1} M + C Q C^\top K^{-1} M - C Q \\
		& + C A \Sigma A^\top C^\top K^{-1} M - C A \Sigma A^\top. 
	\end{align*}
	It shows that the Hessian matrix of $\tr[ M M^\top K^{-1}]$ is positive semi-definite.
\end{proof}

\end{document}